\documentclass{amsart} 

\usepackage{amssymb}

\usepackage{graphicx}

\usepackage{amsmath, amsthm}

\newcommand{\bR}{\mathbb{R}}
\newcommand{\sP}{\mathsf{P}}
\newcommand{\bP}{\mathbb{P}}
\newcommand{\bN}{\mathbb{N}}
\newcommand{\bL}{\mathbb{L}}
\newcommand{\cC}{\mathcal{C}} 
\newcommand{\cS}{\mathcal{S}}
\newcommand{\cO}{\mathcal{O}}
\newcommand{\diam}{\operatorname{diam}}
\newcommand{\sH}{\mathsf{H}}
\newcommand{\Vv}{\boldsymbol{v}}
\newcommand{\Vw}{\boldsymbol{w}}
\newcommand{\bC}{\mathbb{C}}
\newcommand{\FS}{\operatorname{FS}}
\newcommand{\rd}{\mathrm{d}}
\newcommand{\cE}{\mathcal{E}}
\newcommand{\supp}{\operatorname{supp}}
\newcommand{\can}{\operatorname{can}}
\newcommand{\ord}{\operatorname{ord}}
\newcommand{\Res}{\operatorname{Res}}

\newcommand{\ordRes}{\operatorname{ordRes}}
\newcommand{\MinResLoc}{\operatorname{MinResLoc}}
\newcommand{\Crucial}{\operatorname{Crucial}}
\newcommand{\SL}{\mathrm{SL}} 
\newcommand{\PGL}{\mathrm{PGL}}
\newcommand{\bZ}{\mathbb{Z}}
\newcommand{\Rat}{\operatorname{Rat}}
\newcommand{\Proj}{\operatorname{Proj}}
\newcommand{\Spec}{\operatorname{Spec}}
\newcommand{\bA}{\mathbb{A}}
\newcommand{\crit}{\operatorname{crit}}
\newcommand{\bQ}{\mathbb{Q}}
\newcommand{\Fix}{\operatorname{Fix}}
\newcommand{\bif}{\operatorname{bif}}
\newcommand{\bD}{\mathbb{D}}
\newcommand{\cA}{\mathcal{A}}
\newcommand{\cM}{\mathcal{M}}
\newcommand{\cD}{\mathcal{D}}
\newcommand{\sA}{\mathsf{A}}

\newcommand{\rM}{\mathrm{M}}
\newcommand{\rL}{\mathrm{L}}
\newcommand{\widevec}[1]{\overrightarrow{#1}}

\theoremstyle{plain}
\newtheorem{theorem}{Theorem}

\newtheorem{mainth}{Theorem}
\newtheorem{maincoro}{Corollary}

\theoremstyle{definition}
\newtheorem{definition}[theorem]{Definition}
\newtheorem{notation}[theorem]{Notation}
\newtheorem*{convention}{Convention}
\newtheorem*{question}{Question}
\newtheorem*{acknowledgement}{Acknowledgement}
\newtheorem*{conjecture}{Conjecture}
\newtheorem{example}{Example}
\theoremstyle{remark}
\newtheorem{remark}{Remark}

\numberwithin{equation}{section}

\begin{document}

\title[non-archimedean dynamics]{A dynamical system over a non-archimedean field}


\author[Y\^usuke Okuyama]{Y\^usuke Okuyama}
\address{
Division of Mathematics,
Kyoto Institute of Technology, Sakyo-ku,
Kyoto 606-8585 Japan}
\email{okuyama@kit.ac.jp} 



\date{\today}



\begin{abstract} 
This is an expository article, originally written in Japanese, on a dynamical system over
a non-archimedean field. The main viewpoint is from complex and non-archimedean
potential theories. After quickly introducing the Berkovich projective line,
the dynamical moduli space as a scheme, and the various height functions on
the space of rational functions and on the dynamical moduli space,
we first survey our study of Rumely's 
new equivariants in non-archimedean dynamics and then survey
our complex geometric and arithmetic studies of the dynamical moduli space 
from our joint works with Thomas Gauthier and Gabriel Vigny. The latter include
a precise version of McMullen's finiteness theorem on formally exact multiplier spectra
and an effective solution of Silverman's conjecture on a comparison between
the moduli height and the critical height (qualitatively, 
the Silverman-Ingram theorem). The final topic is a degeneration
of complex dynamics.
\end{abstract}

\maketitle


This is an expository article on a dynamical system over
a non-archimedean field, so it would be a good idea to 
begin with the definition of a non-archimedean 
(commutative) field.

\begin{definition}
An absolute value or a norm $|\cdot|$ is said to be
non-archimedean if the strong triangle inequality
\begin{gather*}
 |a+b|\le\max\{|a|,|b|\}
\end{gather*}
holds for any $a,b\in k$. Then the normed (or valued) field $k$
(or $(k,|\cdot|)$) is called a non-archimedean field
(see, e.g.\, the book \cite[\S IV]{NagatafieldEng}).
\end{definition}
The $p$-adic field $\bQ_p$ with 
a $p$-adic norm (from arithmetic) and the field
$\bC((t))$ of Laurent series with a $t$-adic norm 
(from, e.g., complex analytic geometry)
are examples of non-archimedean fields. 
Any field $k$ with a trivial absolute value 
$|\cdot|$ (i.e., $\{a\in k:|a|=1\}=k\setminus\{0\}$)
is also a very useful non-archimedean field.

In this notes\footnote{There are many things to write.}, 
we would like to glance over a few topics 
from a dynamical system induced by a rational function 
(in one variable) of degree $>1$ on the (Berkovich)
projective line defined over a non-archimedean field.
The main viewpoint will be from complex and non-archimedean
potential theories.

\tableofcontents

\vspace*{-20pt}

This expository article is organized as follows.
The first three sections are quick introductions of
the Berkovich projective line, and our main references are
the book \cite{BR10}\footnote{The more recent book \cite{BenedettoBook}
is also recommended.}, 
the survey \cite{Jonsson15}, and the paper \cite{FR09}. 
Section \ref{sec:moduli} also includes a quick introduction
of the dynamical moduli space of complex dynamics
as a scheme, following the paper \cite{Silverman98}, and
a quick introduction of various height functions on the space
of rational functions and on the dynamical moduli space.
In Section \ref{sec:tropical}, 
from a Berkovich hyperbolic geometric viewpoint,
we survey our study \cite{Okugeometric} of Rumely's 
new equivariants in non-archimedean dynamics. 
In Sections \ref{sec:GOV1} and \ref{sec:GOV2},
after quickly introducing
the bifurcation theory in the dynamical moduli space,
we survey our complex geometric and arithmetic studies of the dynamical
moduli space from joint works \cite{GOV16,GOV18}.
The theme in Section \ref{sec:degenerate} is a degeneration
of complex dynamics.

\begin{convention}
 For familiarity, 
 we assume that a non-archimedean field denoted by $K$
 (or $(K,|\cdot|)$) is algebraically 
 closed, that the absolute value 
 $|\cdot|$ of $K$ is non-trivial (i.e., not trivial), 
 and that $K$ is complete as a normed space\footnote{Extending a valued field to such a field if necessary. An algebraically closed and complete field 
 equipped with a non-trivial norm which is not non-archimedean
 is isomorphic to $\bC$ (and the converse is clear).}.
 On the other hand, the characteristic of $K$ 
 is arbitrary unless we mention about that.
\end{convention}

\section{Berkovich projective line and its upper half space}

In this section, we recall 
the definition and basic structures
on a Berkovich projective line.

\subsection{Topological and ordered structures}

Let $\cA$ be a commutative ring having $1$.
A non-negative $\bR$-valued function $[\,\cdot\,]$ on $\cA$ is called a 
(submultiplicative) seminorm on $\cA$ if $[0]=0,[1]=1$, 
and $[\phi\psi]\le[\phi][\psi]$ and $[\phi+\psi]\le[\phi]+[\psi]$ for any $\phi,\psi\in\cA$, and is called a norm on $\cA$
if in addition $\{\phi\in\cA:[\phi]=0\}=\{0\}$.
A seminorm $[\,\cdot\,]$ on $\cA$ is said to be multiplicative if $[1]=1$ and, for any $\phi,\psi\in\cA$,
$[\phi\psi]=[\phi][\psi]$. 
Suppose now that $\cA$ is a normed space equipped with a 
norm $\|\cdot\|$. A seminorm $[\,\cdot\,]$ on $\cA$ is said to be bounded if there is $C\ge 0$ such that 
$[\,\cdot\,]\le C\|\cdot\|$ on $\cA$; 
when $[\,\cdot\,]$ is multiplicative,
we can choose $1$ as the constant $C\ge 0$. 
Suppose now that
$\cA$ is a Banach ring, i.e., 
the normed space $\cA$ (or $(\cA,\|\cdot\|)$)
is a Banach space. 
The set $\cM(\cA)$ of all bounded multiplicative seminorms on $\cA$ is called the spectrum of $\cA$,
which is equipped with such a weakest topology that
for every $\phi\in\cA$, 
the non-negative $\bR$-valued function 
$[\,\cdot\,]\mapsto[\phi]$ is continuous on $\cM(\cA)$.
It is known that $\cM(\cA)$ is non-empty
and is a compact Hausdorff topological space.
This spectrum $\cM(\cA)$ of $\cA$
is also equipped with a (partial) ordering $\le$
(which means either $<$ or $=$) 
so that for any $[\,\cdot\,]_1,[\,\cdot\,]_2\in\cM(\cA)$,
we say $[\,\cdot\,]_1\le[\,\cdot\,]_2$ 
if $[\phi]_1\le[\phi]_2$ for every $\phi\in\cM(\cA)$
(see, e.g., the books \cite[\S1.1, \S1.2]{Berkovichbook}
and \cite[\S C.2]{BR10}).

Now fix a non-archimedean field $K$ 
(or $(K,|\cdot|)$; 
recall our convention on $K$ at the beginning).
A subset in $K$ is called a $K$-closed disk if
$K$ is written as
\begin{gather*}
 B(a,r):=\{z\in K:|z-a|\le r\}\subset K
\end{gather*}
for some $a\in K$ and some $r\ge 0$ (notice that
$B(a,r)=\{a\}$ if (and only if) $r=0$);
here and below, $\subset$ means $\subseteq$. 
By the strong triangle inequality for the $K$,
\begin{itemize}
 \item $B(a,r)=B(b,r)$ for every $b\in B(a,r)$, namely, 
every $b\in B(a,r)$ is a ``center'' of $B(a,r)$,
 \item the ``radius'' $r$ of $B(a,r)$ equals
the (genuine) diameter $\diam B(a,r):=\sup\{|x-y|:x,y\in B(a,r)\}$
of $B(a,r)$ in $(K,|\cdot|)$, and
 \item for any $K$-closed disks $B,B'$, we have
       either $B\subset B'$ or $B\supset B'$ if $B\cap B'\neq\emptyset$.
\end{itemize}

The topological space $K$ 
is totally disconnected. 
This topological issue on $K$ is resolved
by introducing the Berkovich affine line 
defined over $K$.

\begin{definition}[Berkovich closed disk]
For every $R>0$, the generalized Tate algebra 
(defined over $K$ and in one indeterminant $T$)
is the convergent power series ring on $B(0,R)$ 
\begin{gather*}
 K\langle R^{-1}T\rangle:=\Bigl\{\phi(T)=\sum_{j=0}^\infty a_jT^j\in K[[T]]:\lim_{j\to+\infty}|a_j|R^j=0\Bigr\}
\end{gather*}
defined over $K$, which is a Banach ring
equipped with the maximum norm
\begin{gather*}
 \|\phi\|_{B(0,R)}:=\max_j(|a_j|R^j)\Bigl(=\sup_{z\in B(0,R)}|\phi(z)|\Bigr);
\end{gather*}
then $\|\cdot\|_{B(0,R)}$ is multiplicative, and when $R=1$, $\|\cdot\|_{B(0,1)}$ 
is called the Gauss norm on 
$K\langle T\rangle$. A Berkovich closed disk 
defined over $K$ is the spectrum 
\begin{gather*}
 \cD(0,R):=\cM(K\langle R^{-1}T\rangle) 
\end{gather*}
of $K\langle R^{-1}T\rangle$ for some $R>0$, and then
$\|\cdot\|_{B(0,R)}$
is the unique maximal element in the ordered set 
$(\cD(0,R),\le)$ (see, e.g.,
the books \cite[\S1.4.4]{Berkovichbook} and \cite[\S1.2]{BR10}).
\end{definition}

Any element $[\,\cdot\,]\in\cD(0,R)$ restricts to the norm $|\cdot|$
on $K$ and is non-archimedean in that
$[\phi+\psi]\le\max\{[\phi],[\psi]\}$ for any 
$\phi,\psi\in K\langle R^{-1}T\rangle$. 
Remarkably, the following holds

\begin{remark}[Berkovich's representation]\label{th:Berkorep}
For every $[\,\cdot\,]\in\cD(0,R)$, there is a non-increasing
and nesting sequence $(B_n)_{n\in\bN}$ of $K$-closed disks 
such that 
\begin{gather*}
 [\phi]=\inf_{n\in\bN}\sup_{z\in B_n}|\phi(z)|,\quad
\phi\in K\langle R^{-1}T\rangle.
\end{gather*}
\end{remark}

\begin{example}
 To each point $a\in B(0,R)$ is associated 
 the evaluation seminorm $\phi\mapsto[\phi]_a:=|\phi(a)|$ 
 on $K\langle R^{-1}T\rangle$,
 which belongs to  $\cD(0,R)$.
 In representing $[\,\cdot\,]_a$ 
 by a sequence of $K$-closed disks $B_n$ as above,
 we can choose $B_n\equiv \{a\}=B(a,0)$. 
\end{example}

In the following, we fix an affine coordinate $z$ of $\bP^1=\bP^1(K)$, and write $\bP^1=\bA^1\cup\{\infty\}$.

\begin{definition}[Berkovich affine line]
The Berkovich affine line $\sA^1=\sA^1(K)$ 
defined over $K$ is the set of all multiplicative 
seminorm on the polynomial ring $K[z]$ 
restricting to the norm $|\cdot|$ on $K$;
each element $\cS\in\sA^1$ is also written as 
$[\,\cdot\,]_{\cS}$ so that $[\phi]_\cS$ is the value
of $\cS=[\,\cdot\,]_{\cS}$ at $\phi\in K[z]$.
The topology
of $\sA^1$ is the weak topology, which is 
such a weakest topology on $\sA^1$ that for every
$\phi\in K[z]$, the $[0,+\infty)$-valued
function $\cS\mapsto[\phi]_{\cS}$ is continuous
on $\sA^1$. 
\end{definition}
Noting that polynomials are regarded as
power series and that the partial sums of power series
are polynomials, we canonically regard as
\begin{gather*}
 \sA^1=\bigcup_{R>0}\cD(0,R),
\end{gather*}
and it turns out that $\sA^1$ is a locally compact,
locally arcwise connected, uniquely arcwise connected
Hausdorff topological space
(see, e.g., the book \cite[\S2.1]{BR10}).

The Berkovich projective line $\sP^1$ is introduced
both as a topological space 
and as a (partially) ordered set.

\begin{definition}[Berkovich projective line]
Each element $[\,\cdot\,]\in\sA^1$ extends to an 
$[0,+\infty]$-valued function on the rational
function field $K(z)$ so that 
$[\phi]=[\phi_1]/[\phi_0]\in[0,+\infty]$
for each $\phi=\phi_1/\phi_0\in K(z)$, where
$\phi_0,\phi_1\in K[z]$ are coprime. 
Corresponding to $\infty\in\bP^1$,
the $[0,+\infty]$-valued function say $[\,\cdot\,]_\infty$
on $K(z)$ is defined as $[\phi]_\infty=|\phi(\infty)|$,
$\phi\in K(z)$, which is still multiplicative
in an appropriate sense and restricts to the norm $|\cdot|$ on $K$
(see, e.g., the survey \cite[\S3.4]{Jonsson15}).

As a topological space,
the Berkovich projective line defined over $K$ is
\begin{gather*}
 \sP^1=\sP^1(K):=\sA^1\cup\{[\,\cdot\,]_\infty\},
\end{gather*}
equipped with the weak topology, i.e.,
such a weakest topology on $\sP^1$ 
that for every $\phi\in K(z)$, 
the $[0,+\infty]$-valued function
$\cS\mapsto[\phi]_{\cS}$ on $\sP^1$ 
is continuous; then $\sP^1$ is identified with
the one-point compactification of $\sA^1$ 
regarding $\infty=[\,\cdot\,]_\infty$ 
as the additional one point. 

The Berkovich projective line $\sP^1$ is also equipped with 
a (partial) ordering $\le_\infty$ so that for any 
$[\,\cdot\,]_1,[\,\cdot\,]_2\in K[z]$, we say $[\,\cdot\,]_1\le_\infty[\,\cdot\,]_2$
if $[\phi]_1\le[\phi]_2$ for every $\phi\in K[z]$;
the point $[\,\cdot\,]_\infty$
is maximal in the ordered set 
$(\sP^1,\le_\infty)$, and
for every $R>0$, the ordering $\le_\infty$ restricts to the ordering $\le$
of $\cD(0,R)$ regarding $\sA^1=\bigcup_{R>0}\cD(0,R)$.
\end{definition} 

The following extends the statement in Remark \ref{th:Berkorep},
and is useful to understand $\sP^1$ better.

\begin{remark}[Extended Berkovich's representation]
A (possibly empty)
family $\cE$ of $K$-closed disks is said to be (maximal and) nested
if 
\begin{itemize}
 \item for any $B,B'\in\cE$, one contains the other,
 \item for every $B\in\cE$ and every $K$-closed disk $B'$,
if $B\subset B'$, then $B'\in\cE$,
 \item for every non-increasing sequence $(B_n)$ in $\cE$,
       if $\bigcap_nB_n\neq\emptyset$, then
       ($\bigcap_nB_n$ is also a $K$-closed disk and)
       $\bigcap_nB_n\in\cE$.
\end{itemize}
Let $\mathfrak{E}=\mathfrak{E}_K$ be the set of all
nesting families of $K$-closed disks, which is equipped
with such a (partial) ordering $\le$ that
$\cE_1\le\cE_2$ if $\cE_1\supset\cE_2$, so that the empty family
$\emptyset$ is the maximal element in $\mathfrak{E}$,
and we adopt the convention that 
$\bigcap_{\emptyset}=K$.

To each family $\cE\in\mathfrak{E}$ of $K$-closed disks
is associated such a point $\cS_{\cE}\in\sP^1$ as
\begin{gather*}
 [\phi]_{\cS_{\cE}}:=\inf_{B\in\cE}\sup_{z\in B}|\phi(z)|,\quad\phi\in K(z),
\end{gather*}
under the convention that $\inf_{\emptyset}=+\infty$,
and this correspondence $\mathfrak{E}\ni\cE\mapsto\cS_{\cE}\in\sP^1$ is an
isomorphism between the ordered sets, so in particular that
the empty family $\emptyset\in\mathfrak{E}$ 
is associated to $[\,\cdot\,]_\infty\in\sP^1$
(see, e.g., the survey \cite[\S3.3]{Jonsson15}).  
\end{remark}

\begin{definition}[Berkovich upper half space for $\sP^1$]
The Berkovich upper half space for $\sP^1=\sP^1(K)$ is
\begin{gather*}
 \sH^1=\sH^1(K):=\sP^1\setminus\bP^1.
\end{gather*}
For every $\cS=\cS_{\cE}\in\sP^1$,
the (affine) diameter of $\cS$ is defined by
\begin{gather*}
 \diam\cS:=\inf_{B\in\cE}\diam B\in[0,+\infty],
\end{gather*}
under the convention that $\inf_{\emptyset}=+\infty$,
and we also set
\begin{gather*}
 B_{\cS}:=\bigcap_{B\in\cE}B,
\end{gather*}
under the convention that $\bigcap_{\emptyset}=K$.
\end{definition}
For every $\cS=\cS_{\cE}\in\sP^1$,
when $B_{\cS}\neq\emptyset$, then
$B_{\cS}$ is either a $K$-closed disk or the whole $K$ and 
we have $\diam\cS=\diam(B_{\cS})$. 
The $K$-closed disk $B_{\cS}$ is a singleton in $K$ if and only if $\diam\cS=0$,
and we have $\cS_{\cE}=[\,\cdot\,]_\infty$ if and only if $\diam\cE=+\infty$. Hence
\begin{gather*}
 \sH^1=\{\cS\in\sP^1:\diam\cS\in(0,+\infty)\}.
\end{gather*}

\begin{remark}[Berkovich's classification]
Elements in $\sP^1$ are classified into one and only one of
the four types I, II, III, and IV; 
the set of all type I points equals $\bP^1$, and
the sets of all type II, III, IV points equal respectively
\begin{gather*}
\sH^1_{\mathrm{II}}:=\{\cS\in\sH^1:\diam\cS\in|K^\times|\},\quad
\sH^1_{\mathrm{III}}:=\{\cS\in\sH^1:\diam\cS\not\in|K^\times|\}),\\
\sH^1_{\mathrm{IV}}:=\sH^1\setminus(\sH^1_{\mathrm{II}}\sqcup\sH^1_{\mathrm{II}})=\{\cS\in\sH^1:B_{\cS}=\emptyset\}.
\end{gather*}
\end{remark}

\subsection{Tree structure of $\sP^1$}
\label{sec:upperhalf}

The (ordered and closed) interval $[\cS,\cS']$ in $\sP^1$
from $\cS$ to $\cS'$
is defined by the subset in $(\sP^1,\le_\infty)$ 
of all points between $\cS$ and $\cS'$ 
if $\cS\le\cS'$, and in general by the union
$[\cS,\cS\wedge_\infty\cS']\cup[\cS\wedge_\infty\cS',\cS']$ 
in $\sP^1$, where $\cS\wedge_\infty\cS'\in\sP^1$ is such a unique
point in $\sP^1$ that 
\begin{gather*}
[\cS,\infty]\cap[\cS',\infty]=[\cS\wedge_\infty\cS',\infty].
\end{gather*}
Set also the (left half open) interval
$(\cS,\cS']:=[\cS,\cS']\setminus\{\cS\}$ from $\cS$ to $\cS'$.
\begin{figure}[h]
\includegraphics{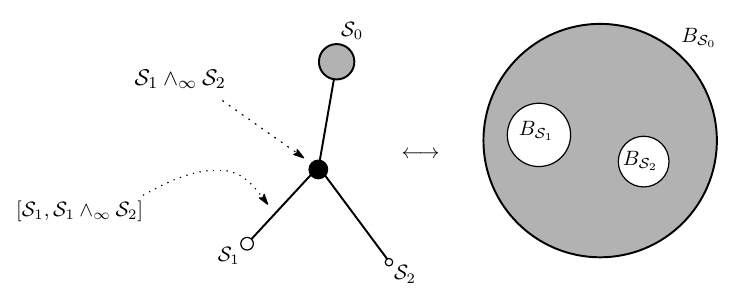}
\caption{$\cS_1\le_\infty\cS_0$ and $\cS_2\le_\infty\cS_0$}
\end{figure}
As an ordered set, the Berkovich $\sP^1$ has the canonical 
(unrooted) ordered tree structure 
in the sense of Jonsson \cite[Definition 2.2]{Jonsson15} 
regarding closed intervals in $\sP^1$ 
as copies of the 
totally ordered subset $[0,1]$ in $(\bR,\le)$ and 
patching them by order preserving isomorphisms.

\begin{definition}[Direction spaces of $\sP^1$ and 
the weak topology of $\sP^1$ as the ordered tree]
For each $\cS\in\sP^1$, the direction space $T_{\cS}\sP^1$ of $\sP^1$ at $\cS$
is the coset of $\sP^1\setminus\{\cS\}$ for such an 
equivalence relation $\sim$ on $\sP^1\setminus\{\cS\}$ that
for any $\cS',\cS''\in\sP^1\setminus\{\cS\}$,
$\cS'\sim\cS''$ if and only if
$(\cS,\cS']\cap(\cS,\cS'']\neq\emptyset$, and 
we denote the equivalence class of $\cS'\in\sP^1\setminus\{\cS\}$ 
by $\widevec{\cS\cS'}$ and call it the
direction from $\cS$ pointing to $\cS'$; the valence 
$v_{\sP^1}(\cS)$ of $\sP^1$ at $\cS$ 
is the cardinality $\#T_{\cS}\sP^1$ 
and takes the value in $\{1,2,+\infty\}$,
and $\widevec{\cS\cS}$ is undefined.
Each direction $\Vv\in T_{\cS}\sP^1$
is also denoted by 
\begin{gather*}
U(\Vv)=\bigl\{\cS'\in\sP^1\setminus\{\cS\}:\widevec{\cS\cS'}=\Vv\bigr\}\end{gather*} 
as a subset
in $\sP^1\setminus\{\cS\}$. 
\end{definition}

The weak topology of $\sP^1$ (as an ordered tree) 
having the quasi-open basis 
$\{U(\Vv):\cS\in\sP^1,\Vv\in T_{\cS}\sP^1\}\subset 2^{\sP^1}$ 
(see, e.g., the survey \cite[\S3.5, \S3.6]{Jonsson15} and
the book \cite[\S2.2]{BR10})
coincides with the (already equipped)
weak topology on $\sP^1$ since a non-empty intersection among a finitely many subsets $U(\Vv)$ in $\sP^1$
for some $\cS\in\sP^1$ and some $\Vv\in T_{\cS}\sP^1$ is nothing but a 
connected open affinoid subset in $\sP^1$
and all of them form an open basis of $\sP^1$. 
Moreover, for every $\cS\in\sP^1$ and every $\Vv\in T_{\cS}\sP^1$,
$U_{\cS}(\Vv)$ is a connected component of $\sP^1\setminus\{\cS\}$.

\begin{example}
We have $U(\widevec{0\infty})=\sP^1\setminus\{0\},
U(\widevec{\infty0})=\sP^1\setminus\{\infty\}=\sA^1$,
so $U(\widevec{0\infty})\cap U(\widevec{\infty0})
=\sA^1\setminus\{0\}$. The point $\cS\in\sP^1$
is of type II, III, IV if and only if
$v_{\sP^1}(\cS)=+\infty,2,1$, respectively.
\end{example}

We see that $\sP^1$ contains both $\bP^1$ and $\sH^1$
as dense subsets, and $\bP^1\cup\sH^1_{\mathrm{IV}}$
is the set of all end points of the tree $\sP^1$.
The weak topology of $\sP^1$ is not necessary 
metrizable.

\begin{definition}[Hyperbolic metric and the Berkovich hyperbolic space]
 The hyperbolic metric on $\sH^1$ is defined by
\begin{gather*}
 \rho(\cS,\cS'):=
\begin{cases}
 \log|\diam(\cS')/\diam\cS|
&\text{if either }\cS\le_\infty\cS'\text{ or }\cS'\le_\infty\cS,\\
\rho(\cS,\cS\wedge_\infty\cS')+\rho(\cS\wedge_\infty\cS',\cS)
&\text{in general};
\end{cases}
\end{gather*}
namely, for any distinct $\cS,\cS'\in\sP^1$, 
$\rho(\cS,\cS')$ is the ``conformal logarithmic modulus''
of the ``Berkovich open annulus'' 
$U(\widevec{\cS\cS'})\cap U(\widevec{\cS'\cS})$.

The metric space $(\sH^1,\rho)$ is called the Berkovich hyperbolic space for $\sP^1$, which is an $\bR$-tree and the 
(ideal) Gromov boundary of which equals $\bP^1$;
the topology of $(\sH^1,\rho)$ is stronger than the
relative one of $\sH^1$ (see, e.g.,
the book \cite[\S2.7]{BR10} and 
the survey \cite[\S3.5]{Jonsson15}).
\end{definition}

We omit the analytic structure of $\sP^1$.
We will see in the following 
that the Berkovich space is useful
in studying various set-theoretic or geometric equivariants in
non-archimedean dynamics, from 
conformal or electrostatic viewpoints as in 
complex dynamics.

\section{Dynamics of rational functions on the Berkovich projective line}

A rational function $h\in K(z)$ defined over a non-archimedean field
$K$ acts on the Berkovich projective line $\sP^1=\sP^1(K)$ so that
for every $\cS\in\sP^1$, the image $h(\cS)\in\sP^1$ is defined so that
\begin{gather*}
 [\phi]_{h(\cS)}=[h^*\phi]_{\cS},\quad \phi\in K(z).
\end{gather*}
This action of $h$ on $\sP^1$ extends
the classical action of $h$ on $\bP^1=\bP^1(K)$, and
is continuous from the definition of the weak topology of $\sP^1$. If $\deg h>0$, 
then the action of $h$ on $\sP^1$ is open and
preserves both $\bP^1=\bP^1$ and 
$\sH^1=\sH^1(K)$, and satisfies
$\#h^{-1}(\cS)\in\{1,\ldots,\deg h\}$ for every $\cS\in\sP^1$; in particular, the action of $h$ on $\sP^1$ is a tree self-map of $\sP^1$ (see \cite[\S2.6]{Jonsson15}).

\subsection{Directional local degree and the local degree}
The action on $\sP^1$ of $h\in K(z)$ of degree $>0$ restricts
to a piecewise affine and $(\deg h)$-Lipschitz 
selfmap of $(\sH^1,\rho)$;
we first note the fact 
that for every $\cS\in\sP^1$ and every $\Vv\in T_{\cS}\sP^1$, 
if $\cS'\in U(\Vv)$ (so that $\Vv=\widevec{\cS\cS'}$)
is close enough to $\cS$, then the action on $\sP^1$ of $h$ 
restricts to an order preserving homeomorphism
$h:[\cS,\cS']\to[h(\cS),h(\cS')]$ between the ordered intervals
and to a proper map
\begin{gather*}
 U(\widevec{\cS\cS'})\cap U(\widevec{\cS'\cS})
\to U(\widevec{h(\cS)h(\cS')})\cap U(\widevec{h(\cS')h(\cS)})
\end{gather*}
between ``Berkovich open annuli''; in particular,
the tangent map $(h_*)_{\cS}:T_{\cS}\sP^1\to T_{h(\cS)}\sP^1$ between the direction spaces
is defined by setting 
\begin{gather*}
 h_*(\Vv):=\widevec{h(\cS)h(\cS')}. 
\end{gather*}
Moreover, when in addition $\cS\in\sH^1$, 
then the restriction 
$h:[\cS,\cS']\to[h(\cS),h(\cS')]$ is even affine with respect to $\rho$-length parametrizations of the ordered intervals
$[\cS,\cS'],[h(\cS),h(\cS')]$ having the slope $m_{\Vv}(h)\in\{1,\ldots,\deg h\}$, so that for any $\cS_1,\cS_2\in[\cS,\cS']$,
\begin{gather*}
 \rho\bigl(h(\cS_1),h(\cS_2)\bigr)=m_{\Vv}(h)\cdot\rho(\cS_1,\cS_2).
\end{gather*}
The integer $m_{\Vv}(h)$ is called the directional local degree 
of $h$ at $\cS$ with respect to $\Vv\in T_{\cS}\sP^1$.
We then note the fact that the local degree function $\deg_{\,\cdot\,}h:\bP^1\to\{1,\ldots,\deg h\}$ of $h$ on $\bP^1$
extends upper semicontinuously to 
the local degree function $\deg_{\,\cdot\,}(h):\sP^1\to\{1,\ldots,\deg h\}$
of $h$ on $\sP^1$ 
so that for every $\cS\in\sH^1$ and every $\Vw\in T_{h(\cS)}\sP^1$,
\begin{gather*}
 \deg_{\cS}(h)=\sum_{\Vv\in T_{\cS}\sP^1:h_*(\Vv)=\Vw}m_{\Vv}(h)
\end{gather*}
and that for every domain (i.e., non-empty connected open subset) $U$
in $\sP^1$ and every component $V$ of $h^{-1}(U)$,
\begin{gather*}
 \cS\mapsto\sum_{\cS'\in h^{-1}(\cS)\cap V}\deg_{\cS'}(h)\equiv
\deg(h:V\to U)\quad\text{on }U,
\end{gather*}
where $\deg(h:V\to U)$ is the degree of the restriction $h:V\to U$
as a proper map (e.g., $\deg(h:\sP^1\to\sP^1)=\deg h$).
From those facts, 
the projective transformation group $\PGL(2,K)$ of $\bP^1$
extends to an isometric automorphism subgroup of $(\sH^1,\rho)$.

The directional local degree and the local degree of $h$
are introduced by Rivera-Letelier \cite[Proposition 3.1]{Juan05} in a geometric way.
Favre and Rivera-Letelier \cite[Proposition-D\'efinition 2.1]{FR09}
introduced the local degree of $h$ using the analytic structure of $\sP^1$, and Baker and Rumely
\cite[\S9]{BR10} introduced both 
the directional local degree and the local degree of $h$
using potential theory (harmonic analysis) on $\sP^1$. 
The directional local degree of $h$
is also introduced 
using the analytic structure of $\sP^1$ by Jonsson \cite[\S 4.6]{Jonsson15}.

\subsection{Equidistribution theorem}\label{sec:equili}

We first recall the equilibrium measure $\mu_f$ and its properties for an endomorphism$f$ of the complex projective space $\bC\bP^k$ of (algebraic)
degree $d>1$ (for the details, see e.g.\ the survey \cite{DSsurvey}).
We denote by $f^n=f^{\circ n}$ the $n$ times iteration of $f$, $n\in\bN$.
From their pluripotential theoretic
study by Fornaess-Sibony \cite{FS94}, there is a weak limit
\begin{gather*}
 \mu_f:=\lim_{n\to+\infty}\frac{(f^n)^*\omega_{\FS}^{\wedge k}}{d^{kn}}
\end{gather*}
on $\bC\bP^k$ (equipped with the Fubini-Study 
K\"ahler form $\omega_{\FS}$ on $\bC\bP^k$
normalized 
so that $\omega_{\FS}^{\wedge k}(\bC\bP^k)=1$).
Namely, letting $\pi:\bC^{k+1}\setminus\{(0,\ldots,0)\}\to\bC\bP^k$ be the canonical projection,
$\|\cdot\|$ the Euclidean norm on $\bC^{k+1}$, and
$\rd\rd^c$ the normalized (complex) Laplacian on $\bC^{k+1}$,
we have $\rd\rd^c\log\|\cdot\|=\pi^*\omega_{\FS}$ on $\bC^{k+1}\setminus\{(0,\ldots,0)\}$.
The probability measure $\mu_f$ on $\bC\bP^k$
is $f$-balanced
in that $f^*\mu_f=d^k\cdot\mu_f$ on $\bC\bP^k$,
has zero masses on any pluripolar subset in $\bC\bP^k$,
and is supported by the $k$-th Julia set $J_k(f)$ of $f$, which is
a subset of the (first) Julia set
\begin{multline*}
 J(f)=J_1(f)\\
:=\bigl\{x\in\bC\bP^k:\text{the iteration family }
(f^n)_{n\in\bN}\text{ is not equicontinuous at }x\bigr\}
\end{multline*}
of $f$ (for the properties of $J(f)$ in the case of $k=1$,
see, e.g., the book \cite{Milnor3rd}). 
Moreover, letting $\cE(f)$ be the maximal
$f$-totally invariant proper algebraic subset in $\bC\bP^k$
\footnote{The existence of $\cE(f)$ is no non-trivial.}
and denoting by $\delta_x$ the Dirac measure on $\bC\bP^k$
at a point $x\in\bC\bP^k$,
the weak convergence $\lim_{n\to+\infty}(f^n)^*\delta_x/d^{kn}$
on $\bC\bP^k$, say the asymptotic equidistribution property for $f$
towards $\mu_f$, holds for any $x\in\bC\bP^k\setminus\cE(f)$.
In particular, $\mu_f$ is mixing under $f$ and is
the unique $f$-balanced probability measure on $\bC\bP^k$ having
no mass on $\cE(f)$, and is indeed the unique maximal entropy measure for $f$ on $\bC\bP^k$ (Fornaess--Sibony \cite{FS94},
Briend--Duval \cite{BD01}, Dinh-Sibony \cite{DS08}). 
Here, the pushforward of a test function $\phi\in C^0(\bC\bP^k)$
is defined by $f_*\phi:=\sum_{y\in f^{-1}(\cdot)}(\deg_y f)\phi(y)\in C^0(\bC\bP^k)$, and the pullback $f^*\nu$ of a measure $\nu$
on $\bC\bP^k$ is defined by the equality
$\int_{\bC\bP^k}\phi(f^*\nu)=\int_{\bC\bP^k}(f_*\phi)\nu$
for every $\phi\in C^0(\bC\bP^k)$.

Coming back to non-archimedean dynamics,
for a rational function $h\in K(z)$ of degree $>0$,
we define the pullback $h^*\nu$ 
of a Radon measure $\nu$ 
on $\sP^1=\sP^1(K)$ in a manner similar to 
that in the last sentence of the previous paragraph; in particular,
for the Dirac measure $\delta_{\cS}$ at a point $\cS\in\sP^1$, we have 
\begin{gather*}
 h^*\delta_{\cS}=\sum_{\cS'\in h^{-1}(\cS)}(\deg_{\cS'}h)\cdot\delta_{\cS'}\quad\text{on }\sP^1,
\end{gather*}
and for a general $\nu$,
we have $h^*\nu=\int_{\sP^1}(h^*\delta_{\cS})\nu(\cS)$.
Let $f\in K(z)$ be of degree $d>1$ and 
write $f^n=f^{\circ n}$ for each $n\in\bN$ as in the previous
paragraph. Corresponding to the statements for 
a morphism of $\bC\bP^k$ of degree $>1$, the following statements
hold; (i) there is a weak limit
\begin{gather*}
 \mu_f:=\lim_{n\to+\infty}\frac{(f^n)^*\delta_{\cS}}{d^n}\quad\text{on }\sP^1
\end{gather*}
for any $\cS\in\sH^1$, (ii)
this probability Radon measure $\mu_f$ on $\sP^1$
is $f$-balanced in that $f^*\mu_f=d\cdot\mu_f$ 
on $\sP^1$
and has no mass on any (potential theoretic)
polar subset in $\sP^1$, (iii) the weak convergence 
$\lim_{n\to+\infty}(f^n)^*\delta_{\cS}/d^n=\mu_f$ on $\sP^1$
holds for every $\cS\in\sP^1\setminus E(f)$,
where $E(f):=\{a\in\bP^1:\#\bigcup_{n\in\bN}f^{-n}(a)<+\infty\}$,
is the (Picard-type) exceptional set of (the iteration family of) $f$
consisting of at most countably many all $f$-totally invariant cycles of $f$,
and (iv) $\mu_f$ is mixing under $f$ and is the unique $f$-balanced
probability Radon measure on $\sP^1$ having no mass on $E(f)$
(for more details, see Baker--Rumely \cite{BR10},
Chambert-Loir \cite{ChambertLoir06},
Favre--Rivera-Letelier \cite{FR09}).

For harmonic analysis on $\sP^1$ used 
in the proof of the above statements, see Section \ref{sec:harmonic}.
A non-archimedean field $K$ treated in this expository article
does not necessarily 
come from arithmetic (so does $\bC_v$
associated to a (finite) place $v$ of a product formula
field $k$). In arithmetic situation,
more precise quantitative statements are obtained
in an electrostatic manner
(see Favre--Rivera-Letelier \cite{FR06}, and
also \cite{OkuDivisor}).

We define the Berkovich Julia set of $f$ as
\begin{gather*}
 J(f):=\supp\mu_f
\end{gather*}
in $\sP^1$. Then $J(f)\cap\bP^1$
coincides with the classical Julia set of 
$(f^n)_{n\in\bN}$ on $\bP^1$, that is,
the locus in $\bP^1$ of non-equicontinuity of $(f^n)_{n\in\bN}$ 
with respect to the chordal metric $[z,w]_{\bP^1}$
(for the chordal metric on $\bP^1$, 
see Section \ref{sec:harmonic}).

\section{Harmonic analysis on $\sP^1$ --- elliptic and hyperbolic geometries}\label{sec:harmonic}

A harmonic analysis on $\sP^1$ is introduced 
from the Berkovich hyperbolic space $(\sH^1,\rho)$
(for the details, see the book \cite{BR10}).

\subsection{Logarithmic chordal kernel on $\sP^1$}

For familiarity, let us herewith introduce
the elliptic or spherical geometry on $\sP^1$,
which is, as we see below,
not necessarily indispensable.
The chordal metric $[z,w]_{\bP^1}$ on 
$\bP^1=\bP^1(K)=\bA^1\cup\{\infty\}$
is written as
\begin{gather*}
 [z,w]_{\bP^1}=\frac{|z-w|}{\max\{1,|z|\}\max\{1,|w|\}},
\quad z,w\in\bA^1,
\end{gather*}
and is normalized so that $[0,\infty]_{\bP^1}=1$
(the notation $[z,w]_{\bP^1}$ is according to the books
Nevanlinna \cite{Nevan70}, Tsuji \cite{Tsuji59}).
The point 
\begin{gather*}
 \cS_{\can}:=\cS_{B(0,1)}\in\sH^1 
\end{gather*}
corresponding to the (constant family of the) $K$-closed unit disk
$B(0,1)$ is called the Gauss (or canonical) point in $\sP^1$. 
For each $\cS_0\in\sP^1$
and any $\cS,\cS'\in\sP^1$, let $\cS\wedge_{\cS_0}\cS'\in\sP^1$
be the triple point among $\cS,\cS',\cS_0$ in that
$[\cS,\cS']\cap[\cS,\cS_0]\cap[\cS',\cS_0]
=\{\cS\wedge_{\cS_0}\cS'\}$ (for example,
the operation $\wedge_\infty$ has already appeared 
in Subsection \ref{sec:upperhalf}); the point
$\cS\wedge_{\cS_{\can}}\cS'$ is also written as 
\begin{gather*}
 \cS\wedge_{\can}\cS' 
\end{gather*}
for simplicity.
Then there is a unique upper semicontinuous and
separately continuous extension $[\cS,\cS']_{\can}$ of $[z,w]_{\bP^1}$ 
to $\sP^1\times\sP^1$ so that
\begin{gather*}
 -\log[\cS,\cS']_{\can}=\rho(\cS_{\can},\cS\wedge_{\can}\cS'),
\quad\cS,\cS'\in\sH^1,
\end{gather*}
the right hand side in which is noting but the
Gromov product on $(\sH^1,\rho)$ with respect to the
point $\cS_{\can}\in\sH^1$ 
(\cite[\S 2.7]{BR10}, \cite[\S 3.4]{FR06}). 
Fixing the second variable
$\cS'\in\sH^1$, the function $\log[\cdot,\cS']_{\can}$
is continuous on $(\sH^1,\rho)$, 
locally constant on $\sH^1\setminus[\cS_{\can},\cS']$,
and affine on $([\cS_{\can},\cS'],\rho)$ having the constant 
slope either $+1$ or $-1$.
The function $[\cS,\cS']_{\can}$ on $\sP^1\times\sP^1$ is called the generalized Hsia kernel function on $\sP^1$
with respect to $\cS_{\can}$ (\cite[\S 4.3]{BR10}), 
but is not a metric on $\sP^1$
noting that $[\cS,\cS]_{\can}=0$ if and only if $\cS\in\bP^1$.

\subsection{Directional derivation on $\sP^1$ and the Laplacian}
\begin{definition}
The convex or connected hull of
a non-empty closed subset $S$ in $\sP^1$ is
\begin{gather*}
 \Gamma_S:=\bigcup_{\cS,\cS'\in S}[\cS,\cS'].
\end{gather*}
A subtree $\Gamma$ in $\sP^1$ is the subset $\Gamma_S$
in $\sP^1$ for some non-empty closed 
$S\subset\sP^1$, and if $\#S<+\infty$, then
this $\Gamma$ is said to be finite.
For a subtree $\Gamma$ in $\sP^1$ and
a point $\cS\in\Gamma$, set
\begin{gather*}
 T_{\cS}\Gamma:=\Bigl\{\Vv\in T_{\cS}\sP^1:\Vv=\widevec{\cS\cS'}\text{ for some }\cS'\in\Gamma\setminus\{\cS\}\Bigr\},
\end{gather*}
and we say $\cS$ is an end point of $\Gamma$ if
$\#T_{\cS}\Gamma\le 1$; the case of $\#T_{\cS}\Gamma=0$ occurs
if and only if $\Gamma$ is trivial, i.e., $\Gamma=\{\cS\}$.  

For subtrees $\Gamma'\subset\Gamma$ in $\sP^1$,
the inclusion from $\Gamma'$ to $\Gamma$
and the retraction from $\Gamma$ to $\Gamma'$
are denoted by $\iota_{\Gamma',\Gamma}$
and $r_{\Gamma,\Gamma'}$, respectively, and then
both 
are continuous with respect to the relative topologies
of $\Gamma',\Gamma$ from $\sP^1$.
\end{definition}

For a point $\cS\in\sH^1$ and $\Vv\in T_{\cS}\sP^1$,
the (distributional)
directional derivation $\rd_{\Vv}$ 
(with respect to
the line element $\rd\rho$) is induced from
the hyperbolic structure of $(\sH^1,\rho)$.
For an affine function $\phi$ on $([\cS_0,\cS],\rho)$,
we compute
\begin{gather*}
 \rd_{\widevec{\cS_0\cS_1}}\phi
=\lim_{\cS\to\cS_0}\frac{\phi(\cS)-\phi(\cS_0)}{\rho(\cS,\cS_0)},
\end{gather*}
which equals the slope of $\phi$ on $[\cS_0,\cS_1]$
at $\cS_0$. Extending $\rho$ to a ``generalized'' metric
on $\sP^1$ which can take the value $+\infty$ appropriately,
as a generalization of the Laplacians on metrized trees, 
the Laplacian $\Delta_\Gamma$ on a subtree $\Gamma\subset\sP^1$ is 
defined so that the domain of $\Delta_\Gamma$ 
is the space $\operatorname{BDV}(\Gamma)$
of functions $\phi$ of bounded derivative variations on $\Gamma$,
that the range of $\Delta_\Gamma$ is the space of Radon measures
on $\Gamma$, and that, for any subtrees $\Gamma'\subset\Gamma$,
the coherence properties
\begin{gather*}
 \Delta_{\Gamma'}\circ(\iota_{\Gamma',\Gamma})^*=
 (r_{\Gamma,\Gamma'})_*\circ\Delta_\Gamma
\quad\text{and}\quad
\Delta_{\Gamma}\circ(r_{\Gamma,\Gamma'})^*=
(\iota_{\Gamma',\Gamma})_*\circ\Delta_\Gamma
\end{gather*}
respectively hold 
on $\operatorname{BDV}(\Gamma)$ and $\operatorname{BDV}(\Gamma')$.
The Radon measure $\Delta_{\Gamma}\phi$ on $\Gamma$
is approximated by the Laplacians
of continuous and piecewise $C^2$ functions $\tilde{\phi}$
with respect to $\rho$ 
on finite subtrees $\Gamma'\subset\Gamma$ in $\sP^1$, 
which are computed as
\begin{gather*}
 \Delta_\Gamma\tilde{\phi}
=(\rd^2\tilde{\phi})\rd\rho+\sum_{\cS\in\Gamma'}\sum_{\Vv\in T_{\cS}(\Gamma')}(\rd_{\Vv}\tilde{\phi})\delta_{\cS}
\quad\text{on }\Gamma'
\end{gather*}
adopting the convention on the sign of the Laplacians from analysis and noting that $\rd^2=(\rd_{\Vw})^2$
is well defined for all but finitely many $\cS\in\Gamma'$ and
that the line element $\rd\rho$ is regarded as the $1$-dimensional Hausdorff measure 
on $\Gamma'$. 

We write $\Delta=\Delta_{\sP^1}$ for simplicity. Then
the logarithmic chordal kernel $\log[\cS,\cS']_{\can}$
is a fundamental solution of the Laplacian $\Delta$ on $\sP^1$ in that
for every $\cS\in\sP^1$,
\begin{gather*}
 \Delta\log[\,\cdot\,,\cS]_{\can}=\delta_{\cS}-\delta_{\cS_{\can}}
\quad\text{on }\sP^1.
\end{gather*}
Moreover, for a rational function $h\in K(z)$ of degree $>0$,
we have ($h^*(\operatorname{BDV}(\sP^1))\subset
\operatorname{BDV}(\sP^1)$ and) the functoriality
\begin{gather*}
 \Delta\circ h^*=h^*\circ\Delta\quad\text{on }
\operatorname{BDV}(\sP^1). 
\end{gather*}
For the introduction of the Laplacian $\Delta$ on 
the (tree) $\sP^1$, we refer to Favre--Jonsson \cite{FJbook}, 
Baker--Rumely \cite[\S 5]{BR10}, and
Favre--Rivera-Letelier \cite[\S 4.1]{FR06}, 
and for a more thorough study on Berkovich curves (which are graphs),
we refer to Thuillier \cite{ThuillierThesis}. See also Jonsson \cite[\S 2.5]{Jonsson15}.
The opposite sign convention of $\Delta$
is adopted in the book \cite{BR10}. 

\section{A tropical function from Berkovich
dynamics related to (potential) good/semistable reductions}
\label{sec:tropical}

We focus on a problem from non-archimedean dynamics
which is at a glance not related to Berkovich spaces.
Let $K$ be a non-archimedean field. The basic notions
from the valuation ring theory are the following.
\begin{definition}
 The unit $K$-closed disk and the unit $K$-``open'' disk
\begin{gather*}
 \cO_K:=B(0,1)\quad\text{and}\quad\mathfrak{m}_K:=\{z\in K:|z|<1\}
\end{gather*} 
are the ring of $K$-integers and the unique maximal ideal
of it. The field $k=k_K:=\cO_K/\mathfrak{m_K}$ is
called the residual field of $K$.
\end{definition}

The canonical projection $\pi_K:K^2\setminus\{(0,0)\}\to\bP^1
 =\bP^1(K)=\bA^1(K)\cup\{\infty\}$ is
 \begin{gather*}
 \pi_K(z_0,z_1)=\begin{cases}
		 \displaystyle \frac{z_1}{z_0} & \text{if }z_0\neq 0,\\
		 \infty &\text{otherwise}.
		\end{cases}
 \end{gather*}

\subsection{Potential good reduction}

The residue class of a $K$-integer $c\in\cO_K$
modulo $\mathfrak{m}_K$ is denoted by $\tilde{c}\in k$;
more generally, 
for a polynomial $P(z)=\sum_{j=0}^{\deg P}c_jz^j\in\cO_K[z]$
of degree $>0$, 
set $\tilde{P}(\zeta):=\sum_{j=0}^{\deg P}\tilde{c_j}\zeta^j\in k[\zeta]$
of degree $\le\deg P$.

The reduction of a point $a\in\bP^1$ modulo $\mathfrak{m}_K$ is
\begin{gather*}
 \tilde{a}:=[\tilde{z_0}:\tilde{z_1}]\in\bP^1(k),
\end{gather*} 
choosing a lift $(z_0,z_1)\in\pi_K^{-1}(z)$
of $a$ which is minimal in that $(z_0,z_1)\in(\cO_K)^2\setminus
(\mathfrak{m}_K)^2$; more generally,

\begin{definition}[a minimal lift of a rational function on $\bP^1$ \cite{KS09}]
For a rational function $f\in K(z)$ of $\deg f>0$,
a (non-degenerate homogeneous) lift of $f$
is an ordered pair 
$F(X,Y)=(F_0(X,Y),F_1(X,Y))\in(K[X,Y]_{\deg f})^2$ of
homogeneous polynomials such that
\begin{gather*}
 f\circ\pi_K=\pi_K\circ F\quad\text{on }K^2\setminus\{(0,0)\},
\end{gather*}
which is unique up to multiplication in 
the unit group $K^\times=K\setminus\{0\}$ of $K$; a lift $F=(F_0,F_1)$ of $f$
is said to be minimal\footnote{according to 
Kawaguchi--Silverman \cite{KS09}.
Although the word ``normalized'' is an option,
it seems a bit vague and there are possible other
``normalizations'' of $F$, e.g., $|\Res F|=1$.}
if 
\begin{gather*}
 \max\{|c|:c\text{ is a coefficient of 
a monomial of }F_0\text{ or }F_1\}=1.
\end{gather*}
The minimal lift of $f$ belongs to $(\cO_K[X,Y]_{\deg f})^2$
and is unique up to multiplication in 
the unit group $\cO_K^\times$ of $\cO_K$.
\end{definition}

The reduction of a rational function $f(z)\in K(z)$ 
of degree $>0$ modulo $\mathfrak{m}_K$ is
\begin{gather*}
 \tilde{f}(\zeta)
:=\frac{\tilde{F_1}(1,\zeta)/\operatorname{GCD}(\tilde{F_0}(1,\zeta),\tilde{F_1}(1,\zeta))}{\tilde{F_0}(1,\zeta)/\operatorname{GCD}(\tilde{F_0}(1,\zeta),\tilde{F_1}(1,\zeta))}\in k(\zeta),
\end{gather*} 
where $F=(F_0,F_1)$ is a minimal lift  of $f$,
so that $\deg(\tilde{f})\in\{0,\ldots,\deg f\}$ and that
$\deg(\tilde{f})>0$ if and only if $f(\cS_{\can})=\cS_{\can}$.

\begin{definition}\label{th:goodred}
We say a rational function $f\in K(z)$ of degree $>0$
has a good reduction (modulo $\mathfrak{m}_K$) if 
$\deg(\tilde{f})=\deg f$, or equivalently,
the homogeneous resultant
\begin{gather*}
 \Res F:=\det\begin{pmatrix}
	    a_0 & \cdots & a_{d-1} & a_d  &        &    \\
	        & \ddots & \vdots  & \vdots & \ddots &    \\
	        &        & a_0     & a_1  & \cdots & a_d\\
	    b_0 & \cdots & b_{d-1} & b_d  &        &    \\
	        & \ddots & \vdots  & \vdots & \ddots &    \\
	        &        & b_0     & b_1  & \cdots & b_d\\ 
	   \end{pmatrix}\in\cO_K
\end{gather*}
of a minimal lift 
$F=(F_0,F_1)=(\sum_{j=0}^{\deg f} a_jX^{d-j}Y^j,
\sum_{\ell=0}^{\deg f} b_\ell X^{d-\ell}Y^\ell)$
of $f$ belongs to $\cO_K^\times$
(for the homogeneous resultant form $\Res$,
see \cite[\S2.4]{SilvermanDynamics}), 
and say $f$ has a potential
good reduction if there is a projective transformation
$h\in\PGL(2,K)$ such that the conjugation
$h\circ f\circ h^{-1}$ of $f$ has a good reduction.
\end{definition}

The problem to determine algorithmically
whether a given non-archimedean 
rational function of degree $>1$ has a potential good
reduction (or not) has been solved by
Bruin--Molnar \cite{BruinMolnar12} 
(over non-archimedean local field) and
finally by Rumely \cite{Rumely13} (over $K$); 
see also Benedetto \cite{Benedetto14}.

\subsection{Rumely's function}

Let us see more details on Rumely's work above.
Pick a rational function $f\in K(z)$ of degree $d>1$, 
and observe that
the function\footnote{There seems no name of this function.}
\begin{multline*}
\ordRes_f:
\PGL(2,K)\ni h\mapsto\\
\ordRes_f(h):=-\log|\Res(\text{a minimal lift of }h^{-1}\circ f\circ h)|\in\bR_{\ge 0}
\end{multline*}
on $\PGL(2,K)$, which is thus induced from $f$ algebraically, takes
the value $0$ (and then $0$ is the minimum of 
this function $\ordRes_f$) if and only if
$f$ has a potential good reduction (Definition \ref{th:goodred}).
We note that the stabilizer subgroup 
in $\PGL(2,K)$ for $\cS_{\can}$ is $\PGL(2,\cO_K)$,
the right action of which on $\PGL(2,K)$ 
is transitive. In particular, the surjection
$\PGL(2,K)\ni h\mapsto h(\cS_{\can})\in
\sH^1_{\mathrm{II}}=\{\cS\in\sH^1:\diam(\cS)\in|K^\times|\}$
descends to the bijection
\begin{gather*}
 \PGL(2,K)/\PGL(2,\cO_K)\ni h\PGL(2,\cO_K)\mapsto h(\cS_{\can})\in\sH^1_{\mathrm{II}},
\end{gather*}
and in turn, the above function
$\ordRes_f$
descends to a function
\begin{gather*}
\sH^1_{\mathrm{II}}\cong\PGL(2,K)/\PGL(2,\cO_K)
\to\bR_{\ge 0}, 
\end{gather*}
which is still denoted by the same notation $\ordRes_f$,
through the above bijection.
Rumely established that 
$\ordRes_f$ extends convexly and properly 
to $(\sH^1,\rho)$, in which $\sH^1_{\mathrm{II}}$ is dense, 
that this extended $\ordRes_f$ on $\sH^1$ is piecewise
affine on closed intervals in $\sH^1$, 
that $\ordRes_f$ always attains its minimum in $\sH^1$,
and that the minimum locus 
\begin{gather*}
 \MinResLoc_f:=\ordRes_f^{-1}\bigl(\min_{\sH^1}\ordRes_f\bigr)
\end{gather*}
of the extended $\ordRes_f$ on $\sH^1$ is a (possibly trivial)
closed segment in $\sH^1$ each end point of which belongs to
$\sH^1_\mathrm{II}$.

By studying intensively the specific non-trivial
finite subtree
\begin{gather*}
 \Gamma_f:=\bigcap_{a\in\bP^1}\Gamma_{\phi^{-1}(a)\cup\{\text{fixed points in }\bP^1\text{ of }f\}}
\end{gather*}
associated with $f$,
Rumely established further properties of $\ordRes_f$
and solved the problem
on algorithmically determining potential good reduction
mentioned in the final paragraph in the previous subsection. 

\subsection{Crucial function}
Rumely did not write down
the extended $\ordRes_f$ on $(\sH^1,\rho)$ explicitly,
and Rumely's analysis of the extended $\ordRes_f$ on $\sH^1$
was based on careful coordinate changes of $\bP^1$ 
(so of $\sP^1$).

Introducing the $f$-crucial\footnote{according
to Rumely's naming the notions related to $\ordRes_f$}
function 
\begin{gather*}
\Crucial_f(\cS)
 :=\frac{\rho(\cS,\cS_{\can})}{2}
 +\frac{\rho(\cS,f(\cS)\wedge_{\can}\cS)
 -\int_{\sP^1}\rho(\cS_{\can},\cS\wedge_{\can}\cdot)\rd(f^*\delta_{\cS_{\can}})}{d-1}
\end{gather*}
on $(\sH^1,\rho)$ (the integral in the right hand side
is indeed a finite sum), which is thus 
a continuous function on $(\sH^1,\rho)$ 
and is defined globally and Berkovich hyperbolic geometrically, 
the author obtained not only an intrinsic and explicit
Berkovich hyperbolic geometric expression
of the function $\ordRes_f$ on $\sH^1$ but also 
a very useful difference (or base change) formula 
of $\Crucial_f$ (so of $\ordRes_f$) as follows.

\begin{mainth}[{\cite[Theorem 1]{Okugeometric}}]\label{th:crucial}
Let $K$ be an algebraically closed field
that is complete with respect to 
a non-trivial and non-archimedean norm, and
pick a rational function $f\in K(z)$ of degree $d>1$.

{\rm (i)} For every $h\in\PGL(2,K)$,
\begin{gather}
 \Crucial_f(h(\cS_{\can}))
=-\frac{1}{2d(d-1)}
\log\frac{|\Res(\text{a minimal lift of }h^{-1}\circ f\circ h)|}{|\Res(\text{a minimal lift of }f)|},
\label{eq:resultant}
\end{gather}
so that
\begin{gather}
 \ordRes_f
=2d(d-1)\cdot\Crucial_f
 -\log|\Res(\text{a minimal lift of }f)|\quad\text{on }\sH^1.
\label{eq:formulaorder}
\end{gather}

{\rm (ii)} Moreover, the difference (or base change) formula
\begin{multline}
\Crucial_f(\cS)-\Crucial_f(\cS_0)\\
=\frac{\rho(\cS,\cS_0)}{2}
+\frac{\rho(\cS,f(\cS)\wedge_{\cS_0}\cS)
-\int_{\sP^1}\rho(\cS_0,\cS\wedge_{\cS_0}\cdot)\rd(f^*\delta_{\cS_0})}{d-1}
\label{eq:conjugate}
\end{multline}
for $\Crucial_f$ holds on $\sH^1\times\sH^1$.
\end{mainth}

The continuous extendability of $\ordRes_f$
from $\sH^1_{\mathrm{II}}$ to $(\sH^1,\rho)$ and
the piecewise affine property of $\ordRes_f$ 
mentioned above immediately follow from \eqref{eq:formulaorder} in
Theorem \ref{th:crucial}.
The other foundational properties of $\ordRes_f$
including the convexity 
and properness of $\ordRes_f$ on $(\sH^1,\rho)$
also follow by (directional) differentiation of $\Crucial_f$;
the difference formula \eqref{eq:conjugate},
which is regarded as a base point change formula,
is designed for this purpose.

For a non-trivial finite subtree $\Gamma\subset\sP^1$,
the valency measure on $\Gamma$ 
is defined by the (signed) Radon measure 
\begin{gather*}
 \nu_{\Gamma}:=(-2)^{-1}\cdot\sum_{\cS\in\Gamma}(\#T_{\cS}\Gamma-2)\cdot(r_{\sP^1,\Gamma})_*\delta_{\cS}\quad\text{on }\Gamma,
\end{gather*}
which is normalized as $\nu_{\Gamma}(\Gamma)=1$ by the Euler
genus formula, and the $\Gamma$-curvature for $\Crucial_f$ is defined 
by the (signed) Radon measure
\begin{gather*}
 \nu_{f,\Gamma}:=\Delta_{\Gamma}(\Crucial_f|\Gamma)+\nu_{\Gamma}
\quad\text{on }\Gamma,
\end{gather*}
so that we still have $\nu_{f,\Gamma}(\Gamma)=1$. 

Among other useful properties of $\Crucial_f$,
we mention the following geometric slope formula of $\Crucial_f$
(so of $\ord_f$).

\begin{mainth}[for the details, {\cite[Theorem 3]{Okugeometric}}]
Under the same assumption as that in Theorem $\ref{th:crucial}$,
for every point $\cS_0\in\Gamma\cap\sH^1$ and
every direction $\Vv\in T_{\cS_0}\Gamma$,
\begin{gather}
 \rd_{\Vv}\Crucial_f=\frac{1}{2}
-\bigl((\iota_{\Gamma,\sP^1})_*\nu_{f,\Gamma}\bigr)\bigl(U(\Vv)\bigr).\label{eq:slope}
\end{gather}
\end{mainth}

Later, Rumely \cite{Rumely14,Rumely17} 
characterized
the minimum locus $\MinResLoc_f$ of $\ordRes_f$
in two ways. The geometric characterization 
asserts the coincidence of
$\MinResLoc_f$ with the barycenter of 
the $\Gamma_f$-curvature $\nu_{f,\Gamma_f}$ for $\Crucial_f$,
and the analytic characterization asserts the coincidence of
$\MinResLoc_f\cap\sH^1_{\mathrm{II}}$ with
the potential GIT-semistable reduction 
locus (modulo $\mathfrak{m}_K$)
for $f$ in $\PGL(2,K)/\PGL(2,\cO_K)(\cong\sH^1_{\mathrm{II}})$.

The geometric invariant theory (GIT) appearing above
also plays a fundamental role in the next Section.

\section{The arithmetic of the dynamical moduli --- from Silverman's conjecture}
\label{sec:moduli}

The monic and centered quadratic polynomial 
family $P_c(z)=z^2+c\in\bZ[c][z]=\bZ[z,c]$ parametrized
by $c\in\bC$ can be iterated as a polynomial in $z$
and its iterations are still belongs to $\bZ[z,c]$. 
Similarly to this parameter $c$ space $\bC$, the dynamical moduli (space) of 
all rational functions of a given degree $>1$ is formulated.
For the details
of schemes over rings, see e.g.\ the book \cite{Hartshorne77}. 

As usual, $\bA_\bZ^N$ and $\bP_\bZ^N$ are
the affine scheme and the projective space
scheme over $\bZ$ of dimension $N$. For each field $k$,
we fix an algebraic closure $\overline{k}$ of $k$.

In the rest of this section, we fix an integer $d>1$. 

\subsection{The dynamical moduli}
The space of rational functions
(on $\bP^1_\bZ$) of degree $d$ is the affine subscheme
\begin{gather*}
 \Rat_d:=\bP^{2d+1}_{\bZ}\setminus V(\rho_d)
=\Spec(\bZ[a,b][\rho_d^{-1}]_{(0)}),
\end{gather*} 
where the ring
$\bZ[a,b]=\bZ[a_0,\ldots,a_d,b_0,\ldots,b_d]$
is graded by the homogeneous degrees and,
writing the indeterminants as $a=(a_0,\ldots,a_d),b=(b_0,\ldots,b_d)$, we set
\begin{gather*}
\rho_d(a,b):=
\det\begin{pmatrix}
	    a_0 & \cdots & a_{d-1} & a_d  &        &    \\
	        & \ddots & \vdots  & \vdots & \ddots &    \\
	        &        & a_0     & a_1  & \cdots & a_d\\
	    b_0 & \cdots & b_{d-1} & b_d  &        &    \\
	        & \ddots & \vdots  & \vdots & \ddots &    \\
	        &        & b_0     & b_1  & \cdots & b_d\\ 
	   \end{pmatrix}\in\bZ[a,b]_{2d};
\end{gather*}
for an algebraically closed field $\Omega$, 
the projective space $\bP^{2d+1}(\Omega)$ is identified with
the set of all the ratios
\begin{gather*}
 [F_a:F_b]:=\biggl[\sum_{j=0}^da_jX^{d-j}Y^j:
\sum_{\ell=0}^db_\ell X^{d-\ell}Y^\ell\biggr]
\end{gather*}
$(a=(a_0,\ldots,a_d),b=(b_0,\ldots,b_d)\in\Omega^{d+1})$
between $F_a,F_b\in\Omega[X,Y]_d$, and we have
$[F_a:F_b]\in\Rat_d(\Omega)$ if and only if
$\rho_d(a,b)\in\Omega^\times=\Omega\setminus\{0\}$, that is,
each point $[F_a:F_b]\in\Rat_d(\Omega)$ 
is identified with the rational function
$f_{[a:b]}=F_b(1,z)/F_a(1,z)\in\Omega(z)$ of degree $d$. 

The special linear algebraic group (scheme)
$\SL_2$ acts on $\bP^{2d+1}_\bZ$ as conjugation;
for an algebraically closed field $\Omega$,
an element
\begin{gather*}
 \phi=\begin{pmatrix}
 \alpha & \beta\\
 \gamma & \delta
\end{pmatrix}\in\SL_2(\Omega)
\end{gather*}
sends each $f=[F_a:F_b]\in\bP^{2d+1}(\Omega)$
to 
\begin{multline*}
 \phi\circ f\circ\phi^{-1}
=[\alpha F_a(\delta X-\beta Y,-\gamma X+\alpha Y)+\beta F_b(\delta X-\beta Y,-\gamma X+\alpha Y):\\
\gamma F_a(\delta X-\beta Y,-\gamma X+\alpha Y)+\delta F_b(\delta X-\beta Y,-\gamma X+\alpha Y)])\in\bP^{2d+1}(\Omega).
\end{multline*}
The space $\Rat_d$ is contained 
in the semistable locus $(\bP^{2d+1}_{\bZ})^{\operatorname{ss}}$ 
of the $\SL_2$-conjugation action on $\bP^{2d+1}$ and is
stabilized by this $\SL_2$-conjugation action, 
and the geometric quotient scheme
\begin{gather*}
 \rM_d:=\Rat_d/\SL_2\cong
\Spec((\bZ[a,b][\rho_d^{-1}]_{(0)})^{\SL_2}),
\end{gather*}
which is as a set the orbit space for $\Rat_d$ under 
the $\SL_2$-conjugation action and is still an integral affine scheme, 
is called the dynamical moduli (space) of rational functions
on $\bP^1_{\bZ}$ of degree $d$ (\cite[Theorem 2.1]{Silverman98});
for the details of the geometric invariant theory
(GIT), see the book \cite{MFK94}.
For an algebraically closed field $\Omega$, as sets,
the affine variety $\rM_d(\Omega)$ coincides with the set 
$\Rat_d(\Omega)/\PGL(2,\Omega)$ of all
$\PGL(2,\Omega)$-conjugacy classes of rational 
functions on $\bP^1(\Omega)$ of degree $d$.

\subsection{Height functions}\label{sec:heights}
Several height functions on $\Rat_d(\overline{k})$ 
and $\rM_d(\overline{k})$ 
are introduced (for height functions,
see, e.g., the books \cite{HindrySilverman00, BombieriGubler06}), 
where $k$ is a product formula field defined as follows. 
Consequently, the spaces $\Rat_d$ and $\rM_d$ posses
significant boundedness and finiteness properties.

\begin{definition}\label{th:product}
 We say a field $k$ is a product formula field
if $k$ is equipped with a set $M_k$
of places $v$ for $k$, i.e., each of which 
is the equivalence class of some non-trivial norm on $k$, with a family
$(|\cdot|_v)_{v\in M_k}$ of representatives of places $v\in M_k$, 
and with a family $(N_v)_{v\in M_k}$
in $\bN$ such that, for every $z\in k^\times=k\setminus\{0\}$,
there is a finite subset $E_z\subset M_k$ 
for which we have $|z|_v=1$ for every $v\in M_k\setminus E_z$,
and that for every $z\in k^\times$, 
the product formula
\begin{gather*}
 \prod_{v\in M_k}|z|_v^{N_v}=1
\end{gather*}
holds. Then a place $v\in M_k$ is said to be finite
if $|\cdot|_v$ is non-archimedean, and otherwise,
said to be infinite.
\end{definition}

According to a theory of valuation fields,
any finite extension $k'$
of a product formula field $k$ is a product formula
field equipped with $M_{k'},(|\cdot|_w)_{w\in M_{k'}}$, and $(N_w)_{w\in M_{k'}}$ 
canonically induced
from the corresponding families for $k$; each $w\in M_{k'}$ is an extension of
some $v\in M_k$, that is written as $w|v$, so that
for each $v\in M_k$, we have not only
$w|v$ for at most finitely
many $w\in M_{k'}$ but also a compatibility property
\begin{gather*}
 \sum_{w|v}N_w=[k':k]
\end{gather*}
on $(N_v)_{v\in M_k}$ and $(N_w)_{w\in M_{k'}}$.

\begin{example}
Equipping $\bQ$ with the equivalence classes $v$
of normalized $p$-adic norms and the Euclidean norm
and the family $(N_v)$ consisting of only $1$,
$\bQ$ becomes a product formula field, so that
any number field (a finite extension of the (product
formula) field $\bQ$) is also a product formula field.

Number field is characterized as a product formula 
field having an infinite place \cite{Artinalgebraic}.
An (algebraic) function field is a product formula field having only finite places.
\end{example} 

Let $k$ be a product formula field. From arithmetic,
the Weil height function on $\bP^N(\overline{k})$ is 
\begin{gather*}
 h_{\bP^N,k}(x):=\frac{\sum_{v\in M_{k'}}N_v\cdot\log\max_{j\in\{0,\ldots,N\}}|x_j|_v}{[k':k]},\quad x=[x_0:\cdots:x_N]\in\bP^N(\overline{k}),
 \end{gather*}
where $k'$ is any finite extension $k'$ of $k$ satisfying $x\in\bP^N(k')$. 
The restriction of $h_{\bP^{2d+1},k}$ 
to $\Rat_d(\overline{k})\subset\bP^{2d+1}(\overline{k})$ defines a height function
\begin{gather*}
 h_{d,k}(f):=h_{\bP^{2d+1},k}([a:b]),\quad f=[a:b]\in\Rat_d(\overline{k})
\end{gather*}
on $\Rat_d(\overline{k})$. 

\begin{definition}[several dynamical height functions]\label{th:heights}
The minimum height function
on $\rM_d(\overline{k})$ is defined by
\begin{gather*}
 h_{d,k}^{\min}([f]):=\min_{g\in[f]}h_{d,k}(f),
\quad[f]\in\rM_d(\overline{k}).
\end{gather*}
Fixing an embedding $\iota:\rM_d\to\bA_{\bZ}^N\subset\bP_{\bZ}^N$, 
we denote by $D$ the restriction of the line bundle
$\cO_{\bP^n(1)}$ on $\bP^N_{\bZ}$
to the Zariski closure say $\overline{\rM_d}$
of $\iota(\rM_d)$ in $\bP^N_{\bZ}$. 
The ample height function $h_{\rM_d,D,k}$ 
on $\rM_d(\overline{k})$ associated to $D$ is
defined by (pulling back by $\iota$) the ample height function
on $\overline{\rM_d}(\overline{k})$ associated to $D$. 

For an individual rational function $f\in\Rat_d(k)$,
the Call--Silverman $f$-height function on $\bP^1(\overline{k})$
is
\begin{gather*}
 \hat{h}_{f,k}(z):=\lim_{n\to+\infty}\frac{h_{\bP^1,k}(f^n(z))}{d^n},\quad z\in\bP^1(\overline{k})
\end{gather*}
(\cite{CS93}, the case of higher dimensional polarized projective varieties is similar). The critical height function
on $\Rat_d(\overline{k})$ is defined by the function
\begin{gather*}
 h_{\crit,k}(f):=\frac{\sum_{c\in C(f)}\hat{h}_{f,k'}(c)}{[k':k]},\quad f\in\Rat_d(\overline{k}),
\end{gather*}
where $k'$ is any finite extension of $k$ 
such that $f\in\Rat_d(k')$ and that the critical set
$C(f):=\{z\in\bP^1(\overline{k}):f'(z)=0\}$
is contained in $\bP^1(k')$, and the sum ranging over $C(f)$ 
in the right hand side in which takes into account the multiplicity of
each $c\in C(f)$ as a critical point of $f$.
The function $h_{\crit,k}(f)$ descends to the critical height function on
$\rM_d(\overline{k})$, which is still denoted by 
$h_{\crit,k}([f])$,
through the projection morphism
$\Rat_d(\overline{k})\to\rM_d(\overline{k})$.
\end{definition}

As we will see in the next section, when $k=\bQ$,
the critical height function $h_{\crit,\bQ}$ on
$\rM_d(\overline{\bQ})$
is a genuine height function on 
$\rM_d(\overline{\bQ})$ except for 
$\rL_d(\overline{\bQ})$, 
where the $d$-th 
flexible Latt\`es locus $\rL_d$ in $\rM_d$
is non-empty if and only if $d=m^2$ for some integer $m>1$,
and then is a (possibly reducible) curve in $\rM_d$
consisting of all $\SL_2$-conjugacy classes
$[f]$ of flexible Latt\`es maps $f\in\Rat_d$.

\subsection{A conjecture by Silverman}\label{th:silvermanconj}

In giving an upper estimates of 
the function $h_{\crit,\bQ}$ on $(\rM_d\setminus\rL_d)(\overline{\bQ})$ below, 
Silverman also conjectured the lower one,
which is later established by Ingram \cite{Ingram18}.

\begin{conjecture}[Silverman {\cite{Silverman12}}]
There are constants $A_1,A_2>0,B_1,B_2\in\bR$ such that
\begin{gather*}
  A_1\cdot h_{\rM_d,D,\bQ}+B_1\le h_{\crit,\bQ}\le A_2\cdot h_{\rM_d,D,\bQ}+B_2\quad\text{on }(\rM_d\setminus\mathrm{L}_d)(\overline{\bQ}).
\end{gather*} 
\end{conjecture}

After Silverman and Ingram,
this conjecture is answered in an effective manner in \cite{GOV18},
pursuing a locally uniform quantitative approximation 
of the Lyapunov exponent $L(f)$
of a rational function $f$ on the projective line
defined over the complex number field $\bC$ or
over a non-archimedean field $K$ (Theorem \ref{th:effectiveapprox});
as a bonus, we obtain an improvement of McMullen's finiteness theorem
(Theorem \ref{th:finiteness}).
Postponing the details of those theorems to the next section,
let us introduce yet another dynamical height
function on $\rM_d(\overline{k})$ in terms of
the (elementary symmetric polynomial functions of the)
multipliers of cycles of rational functions on $\bP^1$.

\subsection{The multiplier height function}\label{sec:multheight}

In terms of the fixed integer $d>1$, for every $n\in\bN$,
we set\footnote{The summation like $\sum_{m\in\bN:m|n}$ is written as
$\sum_{m|n}$ in short when there would be no confusion.}
\begin{gather*}
 d_n:=\sum_{m|n}\mu\Bigl(\frac{n}{m}\Bigr)(d^m+1),
\end{gather*}
where the arithmetic function $\mu:\bN\to\{0,\pm 1\}$
is the M\"obius function (see, e.g., the book 
\cite{Apostol}). We write $\bP^1_{\bZ}=\Proj(\bZ[X,Y])$.
For every $n\in\bN$ and every $f=[F_a:F_b]\in\Rat_d$, 
the $n$-th dynatomic polynomial is defined by
\begin{gather*}
\Phi_n^*(f;X,Y)=\Phi_n^*(a,b;X,Y)
:=\prod_{m|n}\bigl(YF_{a^{(m)}}-XF_{b^{(m)}}\bigr)^{\mu(\frac{n}{m})}\in\bZ[a,b][X,Y]_{d_n},
\end{gather*}
where $f^m=:[F_{a^{(m)}}(X,Y):F_{b^{(m)}}(X,Y)]\in\Rat_{d^m}$.
In the fibered product 
$\bP^1_{\Rat_d}=\bP^1\times_{\Spec\bZ}\Rat_d$,
the zeros of $\Phi_n^*$ defines 
the formally exact period $n$ subscheme
\begin{gather*}
 \Fix_n^{**}:=V(\Phi_n^*), 
\end{gather*}
and similarly, the subscheme
$\Fix_n=V(YF_{a^{(n)}}-XF_{b^{(n)}})$ in $\bP^1_{\Rat_d}$ 
is called the period $n$ subscheme. The projection
$\bP^1_{\Rat_d}\to\Rat_d$ restricts to 
a finite flat morphism $\Fix_n^{**}\to\Rat_d$ of constant degree $d_n$,
and similarly, to a finite flat morphism $\Fix_n\to\Rat_d$ 
of constant degree $d^n+1$
(Silverman \cite[Theorem 4.4]{Silverman98}).
For an algebraically closed field $\Omega$,
as sets, the fiber 
\begin{gather*}
 \Fix^{**}(f^n):=\bigl(\Fix_n^{**}(\Omega)\bigr)_f 
\end{gather*}
of the projection 
$\Fix_n^{**}(\Omega)\to\Rat_d(\Omega)$ 
over $f\in\Rat_d(\Omega)$ is
\begin{multline*}
\bigl\{z_0\in\bP^1(\Omega):f^n(z_0)=z_0,\text{ but } 
f^m(z_0)\neq z_0\text{ for any }m|n, m<n\bigr\}\\
\cup\Bigl\{z_0\in\bP^1(\Omega):\text{ for some } m<n, 
f^m(z_0)=z_0\text{ and}\\
(f^m)'(z_0)\text{ is a primitive }
\frac{n}{m}\text{-th root of unity}\Bigr\}
\end{multline*}
of all periodic points of $f$ in $\bP^1(\Omega)$
having the formally exact period $n$, so that
the union $\bigcup_{m|n}\Fix^{**}(f^m)$ 
is the set of all fixed points of $f^n$ in $\bP^1(\Omega)$.

\begin{notation}
 The relative differential sheaf on a scheme $X$ over 
 a scheme $Y$ is denoted by $\Omega_{X/Y}$.
\end{notation}  
For every $n\in\bN$, $f_{\operatorname{univ}}^n$ is 
the $n$-th iteration of the universal endomorphism 
$f_{\operatorname{univ}}$ of $\bP^1_{\Rat_d}$, 
so that for an algebraically closed field $\Omega$,
the action of $f_{\operatorname{univ}}^n$ on $\bP^1(\Omega)\times\Rat_d(\Omega)$
is a selfmap $(z,f)\mapsto(f^n(z),f)$.
The pullback endomorphism $(f_{\operatorname{univ}}^n)^*$ 
of $\Omega_{\bP^1_{\Rat_d}/\Rat_d}$ restricts
to an $\cO_{\Fix^{**}_n}$-linear endomorphism
on $\Omega_{\bP^1_{\Rat_d}/\Rat_d}|\Fix_n^{**}$,
which we regard as an $\cO_{\Rat_d}$-linear (and diagonal) endomorphism 
of a (locally) free module on $\Rat_d$ of rank $d_n$ 
(recalling that the projection $\Fix_n^{**}\to\Rat_d$ is
finite flat of constant degree $d_n$).
We obtain the eigenpolynomial 
\begin{gather*}
 \deg\bigl(TI_{d_n}-(f_{\operatorname{univ}}^n)^*\bigr)
=T^{d_n}+\sum_{j=1}^{d_n}(-1)^j\sigma_{j,n}^{**}\cdot T^{d_n-j}
\in\bigl(\bZ[a,b][\rho_d^{-1}]_{(0)}\bigr)[T],
\end{gather*}
the coefficients $\sigma_{j,n}^{**}$ in which indeed belong
to $(\bZ[a,b][\rho_d^{-1}]_{(0)})^{\SL_2}$ so descend to
\begin{gather*}
 \sigma_{j,n}^{**}\in H^0(\rM_d,\cO_{\rM_d})\quad (j\in\{1,\ldots,d_n\})
\end{gather*}
through the projection $\Rat_d\to\rM_d$
(Silverman \cite[Theorem 4.5]{Silverman98}).

For a (not necessarily algebraically closed) field $k$, 
indexing $\Fix^{**}(f^n)$ for each $f\in\Rat_d(\overline{k})$
and each $n\in\bN$ as
\begin{gather*}
 z_1(f^n),\ldots,z_{d_n}(f^n)\in\bP^1(\overline{k}) 
\end{gather*}
taking into account the multiplicities of 
each element in $\Fix^{**}(f^n)$ as the
zeros of $\Phi_n^{**}(f;X,Y)$, for each $j\in\{1,\ldots,d_n\}$,
the value $\sigma_{j,n}^{**}(f)$ at $f$ of the above function $\sigma_{j,n}^{**}$ 
on $\Rat_d(\overline{k})$ 
is nothing but the $j$-th elementary
symmetric polynomial of the multipliers
\begin{gather*}
 (f^n)'\bigl(z_1(f^n)\bigr),\ldots,(f^n)'\bigl(z_{d_n}(f^n)\bigr)\in\overline{k} 
\end{gather*}
of $f^n$ at the fixed points $z_1(f^n),\ldots,z_{d_n}(f^n)$ of $f^n$,
and then in fact $\sigma_{j,n}^{**}(f)\in k$.

\begin{definition}[the multiplier morphisms]\label{th:milnor}
For every $n\in\bN$,
the $n$-th formally exact multiplier (or generalized Milnor) morphism 
and the $n$-th multiplier (or generalized Milnor) morphism
on $\Rat_d$ are defined by
\begin{align*}
 s_n^{**}:=&\bigl(\sigma_{1,n}^{**},\sigma_{2,n}^{**},\ldots,\sigma_{d_n,n}^{**}\bigr):\Rat_d\to\bA^{d_n}_{\bZ}\subset\bP^{d_n}_{\bZ},\\
 s_n:=&\bigl(s_m^{**}\bigr)_{m|n}:\Rat_d\to\prod_{m|n}\bA_{\bZ}^{d_m}=\bA^{d^n+1}_{\bZ}\subset\bP^{d^n+1}_{\bZ},
\end{align*}
respectively, both of which descend to $\rM_d$
through the projection $\Rat_d\to \rM_d$.
\end{definition}

\begin{remark}\label{th:milnorsilverman}
 In the case of $d=2$ and $n=1$, we note that
 $d_1=2^1+1=3,2d+1=5$, and $2d-2=2$. 
 A morphism $\rM_2\to\bA^2_{\bZ}$, e.g.,
\begin{gather*}
  (\sigma_{1,1}^{**},\sigma_{1,2}^{**}):\rM_2\to\bA^2_{\bZ}(\subset\bP^2_{\bZ})
\end{gather*}
 obtained by forgetting one of the three components of 
 $s_1^{**}=s_1:\rM_2\to\bA^3_{\bZ}$
 is isomorphic, and extends to an isomorphism
 $\overline{\rM_2}^{\mathrm{s}}\to\bP^2_{\bZ}$.
Here $\overline{\rM_2}^{\mathrm{s}}$
denotes the geometric quotient of
the GIT-stable locus in $\bP^{5}_{\bZ}$ under
the $\SL_2$-conjugation action, which contains $\Rat_2$
(Milnor \cite{Milnor93} over $\bC$, 
Silverman \cite{Silverman98} over $\bZ$).
\end{remark}

As we will see in Section \ref{sec:GOV2}, when $k=\bQ$,
the following yet another dynamical
height functions on $\rM_d(\overline{\bQ})$
are also genuine height functions on 
$(\rM_d\setminus\rL_d)(\overline{\bQ})$.

\begin{definition}[the multiplier height functions]\label{th:multheight}
Let $k$ be a product formula field.
For $n\in\bN$, the $n$-th formally exact multiplier
height function and the $n$-th multiplier height function
are defined by
\begin{gather*}
 \frac{h_{\bP^{d_n},k}\circ s_n^{**}}{n\cdot d_n}\quad
\text{and}\quad
  \frac{h_{\bP^{d^n+1},k}\circ s_n}{n(d^n+1)},
\end{gather*}
respectively, on both $\Rat_d(\overline{k})$ and $\rM_d(\overline{k})$.
\end{definition}

\section{Potential geometry and bifurcation in the dynamical moduli
--- quantitative equidistribution, counting, 
and the volume}
\label{sec:GOV1}

Fix an integer $d>1$. Over $\bC$, the universal endomorphism $f_{\operatorname{univ}}$ of $\bP^1_{\Rat_d}=\bP^1_{\bZ}\times_{\Spec\bZ}\Rat_d$ is regarded as 
the holomorphic family of all rational
functions on $\bP^1(\bC)$ parametrized by 
the complex manifold $\Rat_d(\bC)$.

\subsection{Bifurcation of dynamical systems and its precision using potential geometry}\label{sec:Lyap}
From a general theory on
holomorphic families of rational functions on
$\bP^1(\bC)$ due to Ma\~n\'e--Sad--Sullivan and
Lyubich \cite{MSS,Lyubich83stability}, the space $\Rat_d(\bC)$
is divided into 
\begin{itemize}
 \item  the $J$-stable locus $S_d$ where
 the mapping
 $\Rat_d(\bC)\ni f\mapsto J(f)\in 2^{\bP^1(\bC)}$
 is continuous with respect to the Hausdorff topology
 on the set of compact sets in $\bP^1(\bC)$, which
is an open subset in $\Rat_d(\bC)$, and 
 \item the $J$-unstable (or bifurcation) locus
 $B_d:=\Rat_d(\bC)\setminus S_d$,
which is a non-empty and nowhere dense 
closed subset in $\Rat_d(\bC)$.
\end{itemize}
Moreover, for every $f\in\Rat_d(\bC)$, there are
an open neighborhood $U$ of $f$ and, 
up to taking an at most finitely sheeted
possibly branched holomorphic covering of $U$, 
$2d-2$ holomorphic mapping $c_1,\ldots,c_{2d-2}:U\to\bP^1(\bC)$ such 
that
for every $g\in U$, 
$c_1(g),\ldots,c_{2d-2}(g)\in\bP^1(\bC)$ 
are all the critical points of $g$ taking into
account of their multiplicities and that
$f\in S_d$ if and only if 
for every $j\in\{1,\ldots,2d-2\}$,
$c_j$ is passive\footnote{according to
McMullen's terminology} at $f$ in that
the family $(U\ni g\mapsto g^n(c_j(g))\in\bP^1(\bC))_{n\in\bN}$ of holomorphic mappings is equicontinuous
at the point $f$ in $U$.

Later, the above foundational studies of dynamical
stability and bifurcation for
holomorphic families of  rational functions on $\bP^1(\bC)$ are made more precise at least qualitatively
in a potential-geometric manner (for the details
on pluripotential theory, see e.g.\ the book
\cite{Demaillybook}).
Let us see some details on that 
(but we omit the topic on
the passivity of critical orbits).
For an individual $f\in\Rat_d(\bC)$, the Lyapunov
exponent $L(f)$
of $f$ with respect to the $f$-equilibrium
measure $\mu_f$ on $\bP^1$ 
(see Subsection \ref{sec:equili}) can be defined as
\begin{gather*}
 L(f):=\int_{\bP^1}\log|f'|\mu_f\in[-\infty,+\infty),
\end{gather*}
where $|f'|$ denotes the operator 
norm of the tangent map
$f'$ on $T\bP^1(\bC)$ of $f$ with respect to
any norm on $T\bP^1(\bC)$. 
The Lyapunov exponent
function $L_d:f\mapsto L(f)$ on $\Rat_d(\bC)$ has
its range $[\log\sqrt{d},+\infty)(\subset(0,+\infty))$ 
and is continuous
and plurisubharmonic (by \cite{Ruelle78,Mane88,DeMarco03}, respectively), and descends to the function
\begin{gather*}
 L_d:[f]\mapsto L(f)
\end{gather*}
on $\rM_d(\bC)$ having the same properties as above,
through the holomorphic projection $\Rat_d(\bC)\to\rM_d(\bC)$ from a complex manifold to a complex orbifold. The positive
closed $(1,1)$-current $T_{\bif,d}:=\rd\rd^cL_d$ on $\rM_d(\bC)$ is supported exactly on the image in $\rM_d(\bC)$
of the bifurcation locus $B_d$ in $\Rat_d(\bC)$ under the above projection $\Rat_d(\bC)\to\rM_d(\bC)$
(DeMarco \cite{DeMarco03}). 
For each $p\in\{1,\ldots,2d-2\}$, the positive closed
$(p,p)$-current $T_{\bif,d}^{\wedge p}$ is called the $p$-th bifurcation current on $\rM_d(\bC)$. 
In particular the positive measure
$\mu_{\bif,d}=T_{\bif,d}^{\wedge(2d-2)}$
is called the bifurcation measure on $\rM_d(\bC)$,
which satisfies 
\begin{gather*}
 \int_{\rM_d(\bC)}\mu_{\bif,d}\in(0,+\infty) 
\end{gather*}
(Bassanelli--Berteloot \cite{BassanelliBerteloot07}). 

\begin{figure}[h]
\centering
 \includegraphics[bb=0 0 332 301, width=0.30\textwidth]{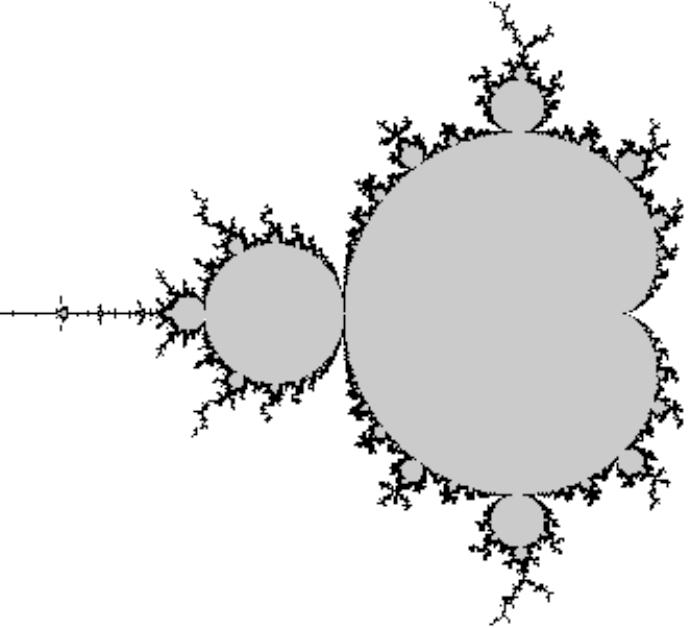}
\caption{The Mandelbrot set $\cC_2$}
\end{figure}
 
Geometrically, 
for each $p\in\{1,\ldots,2d-2\}$,
the support of the $p$-th bifurcation current
on $\rM_d(\bC)$ looks 
locally and generically like a slight distortion of the product of 
$p$ copies of the boundary $\partial\cC_2$ of the Mandelbrot set 
\begin{gather*}
 \cC_2:=\Bigl\{c\in\bC:\limsup_{n\to+\infty}|P_c^n(0)|<+\infty\Bigr\} 
\end{gather*}
(recall that the point $z=0$ is the unique critical point
of $P_c(z)=z^2+c$ in $\bC$ for any parameter $c\in\bC$), 
generalizing McMullen's universality \cite{McMullen2000}
of the (single copy of) $\partial\cC_2$ (Gauthier \cite{Gauthier14}).
This together with Shishikura's equality
$\dim_H(\partial\cC_2)=2(=\dim_\bR\bC)$ recovers 
Gauthier's former result \cite{Gauthier12}, which asserts that
the Hausdorff dimensions of the supports of $T_{\bif,d}^{\wedge p}$ attain the maximal
$4d-4(=\dim_\bR(\bC^{2d-2}))$\footnote{More strongly, the $4d-4$ 
dimensional Lebesgue measure of the support of $\mu_{\bif,d}$ does
not vanish \cite{AGMV19}.}.

\begin{remark}
The bifurcation current (measure)
$\rd\rd^c L(P_\cdot)$ on the parameter $c$ space $\bC$ for the monic and centered 
quadratic polynomial family $(P_c(z))_{c\in\bC}$ coincides with
the harmonic measure  with pole $\infty$ on the above Mandelbrot set 
$\cC_2$.
\end{remark}

\subsection{Locally uniform quantitative approximation of the Lyapunov exponent function} 

In contrast to the compactness of the boundary of the Mandelbrot set $\cC_2$
in $\bC$, the support of the bifurcation measure
$\mu_{\bif,d}$ is not compact in $\rM_d(\bC)$.
Nevertheless, the Lyapunov exponent function
$L_d$ is computed asymptotically
in a locally uniform and quantitative manner as follows;
recall the definitions of $\Fix^{**}(f^n)$ and $\Phi_n^{**}(f;X,Y)$
for each $f\in\Rat_d(\bC)$ and each $n\in\bN$
in Subsection \ref{sec:multheight}.

\begin{mainth}[Gauthier--Okuyama--Vigny \cite{GOV16}]\label{th:approxtruncate}
For a fixed integer $d>1$, there is a constant
$A=A_d\ge 0$ such that
for every $r\in(0,1]$, 
every $f\in\Rat_d(\bC)$, and every $n\in\bN$,
\begin{multline}
\biggl|L(f)-\frac{1}{n\cdot d_n}\sum_{z\in\Fix^{**}(f^n)}\log\max\{r,|(f^n)'(z)|\}\biggr|\\
 \leq A\cdot\Bigl(\log\sup_{\bP^1}(f^\#)+\sup_{\bP^1}|g_f|+|\log r|\Bigr)\frac{\sum_{m|n}m^2}{d^n}.
\label{eq:locunif}
\end{multline}
Here the sum ranging over the
set $\Fix^{**}(f^n)$ of all periodic points $z$
in $\bP^1(\bC)$ of $f$ having the formally exact period $n$
takes into account the multiplicity of each $z$ as a zero of $\Phi_n^{**}(f;X,Y)$.
\end{mainth}
The estimate
\eqref{eq:locunif} is regarded as a 
locally uniform quantitative approximation formula
of the Lyapunov exponent function $L_d$ on $\rM_d(\bC)$ 
in terms of the truncated multipliers of
periodic points of $f$ (or of the $\PGL(2,\bC)$-conjugacy class $[f]$). Here the function
\begin{gather*}
 f^\#(z):=\lim_{w\to z}\frac{[f(z),f(w)]_{\bP^1}}{[z,w]_{\bP^1}},\quad z\in\bP^1(\bC), 
\end{gather*}
with respect to
the chordal metric $[z,w]_{\bP^1}$ on $\bP^1(\bC)$
is the chordal derivative function for $f$
and that $g_f=G^F-\log\|\cdot\|$ on $\bP^1(\bC)$
is the $f$-dynamical Green function on $\bP^1(\bC)$,
where $\|\cdot\|$ is the Euclidean norm on $\bC^2$
and $G^F:=\lim_{n\to+\infty}(\log\|F^n\|)/d^n$
is the escaping rate function for a lift\footnote{Those $F$ and $\Res F$ are defined in a manner similar to
that for non-archimedean fields $K$.} 
$F\in(\bC[X,Y]_d)^2$ of $f$ satisfying $|\Res F|=1$.

In both complex and non-archimedean dynamics,
without truncation (i.e., letting $r=0$),
an approximation formula of $L(f)$ similar to
\eqref{eq:locunif} for an individual $f$ 
has been known with a better error estimate 
$O((\sum_{m|n}m)/d^n)$ as $n\to+\infty$
(\cite{OkuLyap}). The first approximation of
$L(f)$ of this kind (with no non-trivial error estimate) is due to Szpiro--Tucker \cite{ST05}, where they worked over a product formula field
and used Roth's theorem from Diophantine 
approximation.

Our pursuit in Theorem \ref{th:approxtruncate}
of the locally uniform quantitative 
approximation of the Lyapunov exponent function 
$L_d$ by the (truncated) multipliers of periodic points of 
rational functions concludes several 
potential geometric properties
of the bifurcation loci in $\rM_d(\bC)$, as follows. 

\subsection{Quantitative equidistribution, counting, and the volume in $\rM_d$}

The following is the analog of the so called centers
(i.e., the parameters $c\in\bC$ for which 
the unique critical point $z=0$ in $\bC$
of $P_c(z)=z^2+c$ is periodic under $P_c$)
of hyperbolic components of the interior of the Mandelbrot set $\cC_2$.

\begin{definition}[the disjoint and postcritically finite hyperbolic loci in $\rM_d(\bC)$]\label{th:center}
For each $(2d-2)$-tuples 
\begin{gather*}
 \underline{n}=(n_1,\ldots,n_{2d-2})\in\bN^{2d-2},
\end{gather*}
the disjoint and PCF (i.e.\ postcritically finite) hyperbolic locus 
$C_{\underline{n}}$ in $\rM_d(\bC)$
of type $\underline{n}$ is the algebraic set
(defined over $\bQ$) in $\rM_d(\bC)$ 
consisting of all $\PGL(2,\bC)$-conjugacy classes
$[f]\in \rM_d(\bC)$ such that, indexing all 
the $2d-2$ critical points in $\bP^1(\bC)$
(taking into account the multiplicities of them)
of a representative $f$ as $c_1,\ldots,c_{2d-2}$
appropriately, we have both
\begin{description}
\item[PCF hyperbolicity] 
for every $j\in\{1,\ldots,2d-2\}$,
$c_j\in\Fix^{**}(f^{n_j})$, and
\item[Disjointness] 
for any distinct $j,k\in\{1,\ldots,2d-2\}$,
$c_j\not\in\bigl\{f^m(c_k):m\in\{1,\ldots,n_k\}\bigr\}$.
\end{description}
\end{definition}

The fact that $\dim C_{\underline{n}}=0$
is seen by a transversality argument using 
infinitesimal deformation (of Kodaira--Spencer type)
of rational functions, and we have an upper bound of the cardinality
$\#C_{\underline{n}}$ using B\`ezout's
theorem\footnote{This is already non-trivial, but
we omit the detail since we would state a more precise
Theorems \ref{th:equidist} and \ref{th:counting}.}.
\begin{figure}[h]
\centering
\hspace*{-20pt}
\includegraphics[bb=0 0 999 999, width=0.3\textwidth]{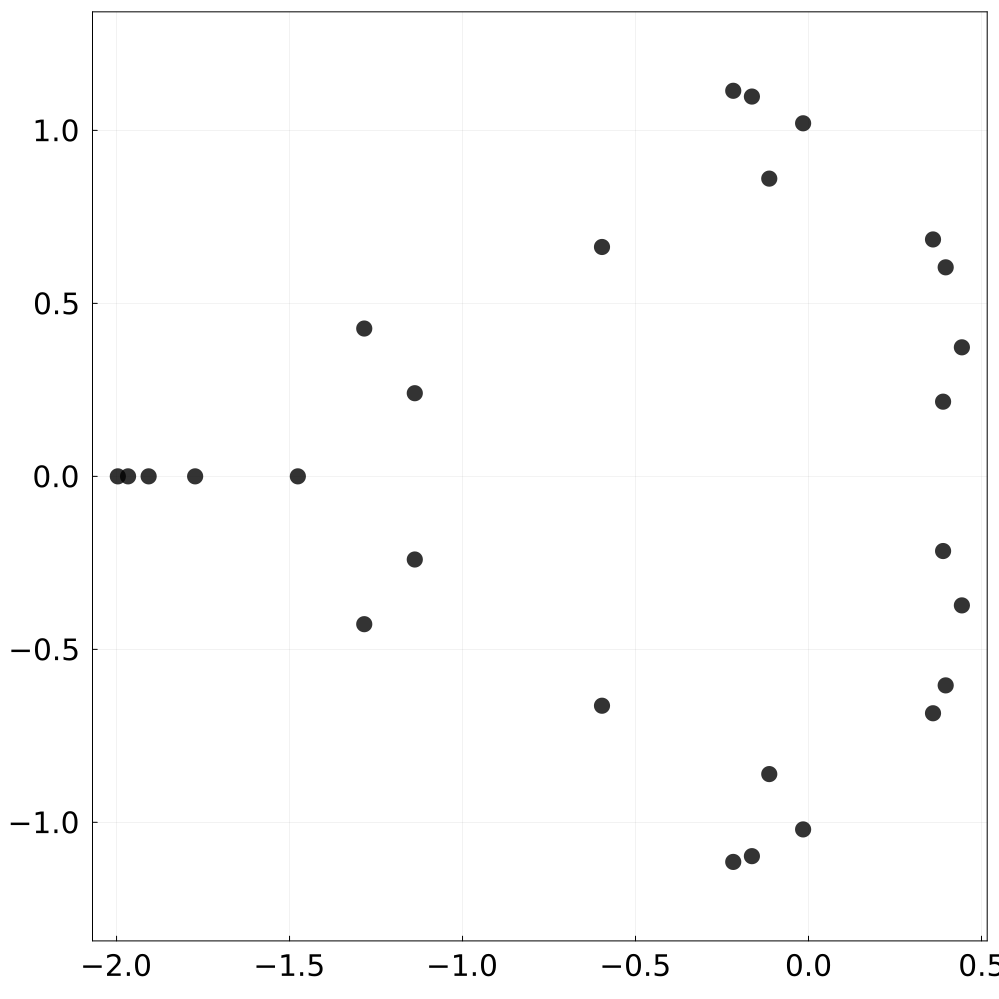}
\caption{The PCF hyperbolic locus $C_6$ in the Mandelbrot set}
\end{figure}

For each $(2d-2)$-tuples $\underline{n}=(n_1,\ldots,n_{2d-2})\in\bN^{2d-2}$, we set
\begin{gather*}
 d_{\underline{n}}:=\prod_{j=1}^{2d-2}d_{n_j}
\end{gather*}
and set
$\operatorname{Stab}(\underline{n})
:=\bigl\{\sigma\in S_{2d-2}:
\bigl(n_{\sigma(1)},\ldots,n_{\sigma(2d-2)}\bigr)=\underline{n}\bigr\}$,
which is a subgroup
of the $(2d-2)$-th symmetric group $S_{2d-2}$.

Using the intersection theory of currents,
from Theorem \ref{th:approxtruncate}, we first deduce
a quantitative equidistribution of
the averaged counting measure
 \begin{gather*}
 \mu_{\underline{n}}:=\frac{\#\operatorname{Stab}(\underline{n})}{d_{\underline{n}}}
 \sum_{[f]\in C_{\underline{n}}}\delta_{[f]}\quad\text{on }\rM_d(\bC)
\end{gather*}
of the disjoint and PCF hyperbolicity locus 
$C_{\underline{n}}$ of type $\underline{n}$ 
towards the bifurcation measure $\mu_{\bif,d}$ on $\rM_d(\bC)$.

\begin{mainth}[Gauthier--Okuyama--Vigny {\cite{GOV16}}, quantitative equidistribution]\label{th:equidist}
Fix an integer $d>1$. 
For every compact subset
$K$ in $\rM_d(\bC)$, there is a constant
$C_K=C_{K,d}>0$ such that for every
$C^2$-test function $\phi$ supported by $K$
and every $(2d-2)$-tuple 
$\underline{n}=(n_1,\ldots,n_{2d-2})\in\bN^{2d-2}$,
\begin{gather*}
 \left|\int_{\rM_d(\bC)}\psi\bigl(\mu_{\underline{n}}-\mu_{\bif,d}\bigr)\right|
\le C_K\cdot\|\psi\|_{C^2}\cdot\max_{j\in\{1,\ldots,2d-2\}}\biggl(\frac{\sum_{m|n_j}m^2}{d_{n_j}}\biggr).
\end{gather*}
In particular,
for any sequence
$(\underline{n}_k)_{k=1}^{\infty}$ in $\bN^{2d-2}$,
we have\footnote{By $\min_j n_{k,j}\to+\infty$, we mean
$\lim_{k\to+\infty}\min_j n_{k,j}=+\infty$ and to make $k\to+\infty$.} 
\begin{gather*}
 \lim_{\min_j n_{k,j}\to+\infty}\mu_{\underline{n}_k}=\mu_{\bif,d}\quad\text{weakly
on }\rM_d(\bC).
\end{gather*}
\end{mainth}

The constant $C_K>0$ in Theorem \ref{th:equidist}
explodes when $K$ gets closer and closer to $(\overline{\rM_d}\setminus\rM_d)(\bC)$,
and the explosion rate is controlled in Theorem \ref{th:approxtruncate}
(or in Theorem \ref{th:effectiveapprox} below more effectively).
By a truncation argument from
pluripotential theory for 
an embedded $\rM_d(\bC)$ in $\bA^N(\bC)\subset\bP^N(\bC)$
and an argument similar to that in the proof of
Theorem \ref{th:equidist}, we also count 
$C_{\underline{n}}$ asymptotically and quantitatively
in terms of the volume $\int_{\rM_d(\bC)}\mu_{\bif,d}$.

\begin{mainth}[Gauthier--Okuyama--Vigny {\cite{GOV16}}, counting and the volume]\label{th:counting}
Fix an integer $d>1$. For every $(2d-2)$-tuple
$\underline{n}=(n_1,\ldots,n_{2d-2})\in\bN^{2d-2}$,
\begin{gather}
 \frac{\#\operatorname{Stab}(\underline{n})}{d_{\underline{n}}}
\cdot \#C_{\underline{n}}=\int_{\rM_d(\bC)}\mu_{\bif,d}
+O\biggl(\max_{j\in\{1,\ldots,2d-2\}}\biggl(\frac{\sum_{m|n_j}m^2}{d_{n_j}}\biggr)\biggr)
\label{eq:countingmeasure}
\end{gather}
as $\min\{n_1,\ldots,n_{2d-2}\}\to+\infty$.
\end{mainth}

When $d=2$, combining the geometric counting of $\cC_{\underline{n}}$ in
\eqref{eq:countingmeasure} with Kiwi--Rees's 
algebraic geometric counting of $C_{\underline{n}}$
\cite{KR13} (based on the isomorphism $\overline{\rM_d}^s(\bC)\cong\bP^2(\bC)$ of Milnor (and Silverman)
mentioned in Remark \ref{th:milnorsilverman}), we 
establish an exact formula of 
the volume of $\rM_2(\bC)$ with respect to $\mu_{\bif,2}$.

\begin{maincoro}[the mass formula of
$\mu_{\bif,2}$ \cite{GOV16}]
Let $\phi$ denote Euler's totient function. Then
\begin{gather}
 \int_{\rM_2(\bC)}\mu_{\bif,2}
=\frac{1}{3}-\frac{1}{8}\sum_{n=1}^{\infty}\frac{\phi(n)}{(2^n-1)^2}.\tag{\ref{eq:countingmeasure}$'$}\label{eq:mass}
\end{gather}
\end{maincoro}

The following might be of some interest.

\begin{question}
 Is the series in the right hand side in
 \eqref{eq:mass} a rational number or not?
\end{question}

\section{Arithmetic of the dynamical moduli --- 
an improvement of McMullen's finiteness theorem, 
and effective comparisons among dynamical height functions}
\label{sec:GOV2}

When $d=2$ and $n=1$, 
an isomorphism $\rM_2\to\bA^2_{\bZ}$
is obtained by forgetting any one component
of the first (formally exact) multiplier morphism
$s_1=s_1^{**}=(\sigma_{1,1}^{**},\sigma_{1,2}^{**},\sigma_{1,3}^{**}):\rM_2\to\bA^3_{\bZ}$, which is in particular injective (Remark \ref{th:milnorsilverman}).

Fix an integer $d>1$.

\subsection{An improvement of McMullen's finiteness theorem}

McMullen's finiteness theorem \cite{McMullen87} asserts that
that there is $N_0\in\bN$ such that 
for every $n>n_0$, any fiber of the restriction
\begin{gather*}
 s_n:(\rM_d\setminus\rL_d)(\bC)\to\bA^{d^n+1}(\bC)
\end{gather*}
of the $n$-th multiplier spectrum 
$s_n=(s_m^{**})_{m|n}:\rM_d(\bC)\to\prod_{m|n}\bA^{d_m}_{\bZ}(\bC)=\bA^{d^n+1}(\bC)$ is finite;
for the $s_n$, the formally exact multiplier spectra $s_m^{**}$, and
$d$-th flexible Latt\`es locus $\rL_d$ in $\rM_d$,
which is non-empty if and only if $d=m^2$ for some 
integer $m>1$, see Subsection \ref{sec:multheight} and
the final paragraph in Subsection \ref{sec:heights}. 
The proof was based on a normal family argument and 
Thurston's rigidity theorem\footnote{For a proof, see Douady--Hubbard \cite{DH93}.}
in complex dynamics.

Recall that 
we denoted 
by $D$ the restriction of 
the line bundle $\cO_{\bP^N}(1)$ to the
Zariski closure $\overline{\rM_d}$ in $\bP^N_{\bZ}$ fixing an 
embedding $\rM_d\to\bA^N_{\bZ}\subset\bP^N_{\bZ}$ 
(Definition \ref{th:heights}). 
Let also $\omega_D$ be the restriction of the Fubini-Study
K\"ahler form on $\bP^N(\bC)$ to $\rM_d(\bC)$.

An improvement of McMullen's finiteness theorem is established
even in an effective manner, developing further the proof of 
our effective version of Silverman-Ingram's comparison theorem 
(see Conjecture in Subsection \ref{th:silvermanconj}), as follows.

\begin{mainth}[Gauthier--Okuyama--Vigny {\cite{GOV18}}]\label{th:finiteness}
Fixing an integer $d>1$,
there is a constant $n_1\in\bN$ such that for every $n\ge n_1$,
any fiber of the restriction 
\begin{gather*}
 s_n^{**}:(\rM_d\setminus\rL_d)(\bC)\to\bA^{d_n}(\bC)
\end{gather*}
of the $n$-th formally exact multiplier morphism
$s_n^{**}=s_{n,d}^{**}:\rM_d(\bC)\to\bA^{d_n}(\bC)$ 
to $(\rM_d\setminus\rL_d)(\bC)$ is $($already$)$
finite. The largeness of the constant $n_1$ is 
effectively determined using only $d$,
the complex analytic quantities 
$\|\mu_{\bif,d}\|_{\rM_d(\bC)}=\int_{\rM_d(\bC)}\mu_{\bif,d}$ and $\|T_{\bif,d}\|_{\rM_d(\bC)}=\int_{\rM_d(\bC)}T_{\bif,d}\wedge\omega_D^{2d-3}$, and the algebraic quantity $\deg_D(\rM_d)$.
\end{mainth}

In the next two sections, let us see the proof outline
of our effective version 
of Silverman-Ingram's comparison theorem.

\subsection{Effective comparison between the multiplier height 
and critical height functions}

For a non-archimedean field $K$, the chordal derivative function
\begin{gather*}
 h^\#(z):=\lim_{\bP^1\ni w\to z}\frac{[h(z),h(w)]_{\bP^1}}{[z,w]_{\bP^1}}:\bP^1=\bP^1(K)\to[0,+\infty)
\end{gather*}
of a rational function $h\in K(z)$ extends 
continuously to $\sP^1=\sP^1(K)$. Fix an integer
$d>1$. For every $f\in\Rat_d(K)$, the Lyapunov 
exponent of $f$ with respect to the $f$-equilibrium 
measure $\mu_f$ on $\sP^1$ (see Subsection \ref{sec:equili}) 
is well-defined by
\begin{gather*}
 L(f):=\int_{\sP^1}\log(f^\#)\mu_f\in\bR;
\end{gather*}
for each $f\in\Rat_d(\bC)$, the Lyapunov exponent
of $f$ with respect to the $f$-equilibrium
measure $\mu_f$ on $\bP^1(\bC)$ (see Subsection \ref{sec:Lyap})
is also written
as $L(f)=\int_{\bP^1(\bC)}\log(f^\#)\mu_f(\in[\log\sqrt{d},+\infty))$ using 
the chordal derivative function $f^\#$ of $f$ on $\bP^1(\bC)$
(see the paragraph after Theorem \ref{th:approxtruncate}).

Let $k$ be a product formula field (see Definition
\ref{th:product}).

\begin{notation}
 For each place $v\in M_k$, the field $\bC_v$ defined by
 the completion (with respect to the extended 
 norm $|\cdot|_v$) of an algebraic closure of
 the completion $k_v$ of $k$ with respect to 
 $|\cdot|_v$ is also algebraically closed,
 and then we fix an embedding of
 the algebraic closure $\overline{k}$ into $\bC_v$.

 Similarly to the notations $|\cdot|_v$, $N_v$, $k_v$, and $\bC_v$, we indicate the dependence of a quantity 
 (mainly induced by $f\in\Rat_d$) on $v$ by adding the suffix $v$ to the notation for
 this quantity; for example, for each $f\in\Rat_d(\overline{k})$ and $v\in M_k$, the quantity $L(f)$ is denoted by $L(f)_v$ when we regard the $f$ as $f\in\Rat_d(\bC_v)$.
\end{notation}

The critical height function on $\rM_d(\overline{k})$
(see Subsection \ref{sec:heights})
is written as
\begin{gather*}
h_{\crit,k}([f])=\frac{1}{[k':k]}\sum_{w\in M_{k'}}N_w\cdot L(f)_w,\quad [f]\in\rM_d(\overline{k})
\end{gather*}
by an integration by parts and the product formula,
where $k'$ is a finite field extension
of $k$ such that $f\in\Rat_d(k')$. On the other hand,
for every $n\in\bN$, decoding 
the (a bit cryptic) definition of the $n$-th formally exact multiplier height function on $\rM_d(\overline{k})$ (see Definition \ref{th:multheight}), we have
\begin{multline*}
 \frac{(h_{\bP^{d_n},k}\circ s_n^{**})([f])}{n\cdot d_n}\\
=\frac{1}{[k':k]}
\sum_{w\in M_{k'}}N_w\cdot\frac{\sum_{z\in\Fix^{**}(f^n)}\log\max\{1,|(f^n)'(z)|_w\}
}{n\cdot d_n},\quad [f]\in\rM_d(\overline{k}),
\end{multline*}
where $k'$ is a finite field extension
of $k$ such that $f\in\Rat_d(k')$ and $\Fix^{**}(f^n)\subset\bP^1(k')$.

The following is a precision of 
Theorem \ref{th:approxtruncate} 
(a locally uniform quantitative approximation formula
of the Lyapunov exponent function $L_d$)
applied not only to $\bC$ but also 
to a non-archimedean field $K$, which is
not necessarily of arithmetic origin, of characteristic $0$.

\begin{mainth}[Gauthier--Okuyama--Vigny \cite{GOV18}]\label{th:effectiveapprox}
Let $K$ be an algebraically closed field
of characteristic $0$ that is complete with
respect to a non-trivial norm $|\cdot|$.
Then fixing an integer $d>1$,
for every $f\in\Rat_d(K)$, 
every $n\in\bN$, and every 
$r\in(0,\epsilon_{d^n}]$,
\begin{multline}
\biggl|L(f)-\frac{1}{n\cdot d_n}\sum_{z\in\Fix^{**}(f^n)}\log\max\{r,|(f^n)'(z)|\}\biggr|\\
 \leq 2(2d-2)^2\left(|L(f)|+16\cdot\frac{3}{2}\log(M_1(f)^2)+\sup_{\sP^1}|g_f|+|\log r|\right)\frac{\sum_{m|n}m^2}{d_n},\label{eq:local}
\end{multline}
and the number $16$ in the right hand side could decrease to $1$
when $K$ is non-archimedean.

Here $\epsilon_d:=\min\bigl\{|m|^d:m\in\{1,\ldots,d\}\bigr\}\in|K^\times|\cap(0,1]$ is
the Benedetto--Ingram--Jones--Levy constant \cite{BIJL14}
for $d$.
\end{mainth}

The constant $M_1(f)$ is defined by
\begin{gather*}
 M_1(f):=\begin{cases}
	  \sup_{\bP^1}(f^\#) & \text{for }K=\bC,\\
	  |\Res(\text{a minimal lift of }f)|
& \text{for non-archimedean }K,
	 \end{cases}
\end{gather*} 
so that $f:\bP^1(K)\to\bP^1(K)$ is $M_1(f)$-Lipschitz
continuous with respect to the chordal metric on $\bP^1(K)$
(due to Rumely--Winburn \cite{RumelyWinburn15}
for non-archimedean $K$).
Similarly to the case $K=\bC$ (see the paragraph after 
Theorem \ref{th:approxtruncate}), for a non-archimedean $K$,
the function $g_f$ is the dynamical 
Green function for $f$
on $\sP^1(K)$ defined by the continuous extension to $\sP^1(K)$ of the function $G^F-\log\|\cdot\|$ on $\bP^1(K)$, where $\|\cdot\|$ is the maximal norm on $K^2$
and $G^F:=\lim_{n\to +\infty}(\log\|F^n\|)/d^n$ on $K^2\setminus\{(0,0)\}$ is the escaping rate function for a lift $F$ of $f$ satisfying
$|\Res F|=1$. The explicit constants in the right hand side
in \eqref{eq:local} (except for $|\log r|$)
are estimated in terms of the algebraic quantities
$|\Res F|$ and $|F|:=\max\{|\text{a coefficient of }F_0\text{ or }F_1|\}$,
where $F=(F_0,F_1)\in(K[X,Y]_d)^2$ is a lift of $f$. 

When the product formula field $k$ is $\bQ$,
summing up \eqref{eq:local} for $K=\bC_v$ over all places $v\in M_{\bQ}$,
we obtain the following effective comparisons
between the critical/(formally exact) multiplier height functions on 
$\rM_d(\overline{\bQ})$, which also involves
the minimal height function on $\rM_d(\overline{\bQ})$;
for every integer $d>1$ and every $n\in\bN$,
there is a constant $C_{d,n}\in\bR$ (depending only on $d,n$) such that 
\begin{multline}
 \left|
\frac{h_{\bP^{d_n},\bQ}\circ s_n^{**}}{n\cdot d_n}
- h_{\crit,\bQ}\right|\\
\leq 8(d-1)(196d^2-192d-3)\cdot h^{\min}_{d,\bQ}\cdot \frac{\sum_{m|n}m^2}{d_n}+C_{d,n}
\quad\text{on }\rM_d(\overline{\bQ}),\tag{\ref{eq:local}$'$}\label{eq:multcritexact}
\end{multline}
and in turn, by M\"obius inversion of \eqref{eq:multcritexact}, 
there is also a constant $C'_{d,n}\in\bR$ such that
\begin{multline}
 \left|\frac{h_{\bP^{d^n+1},\bQ}\circ s_n}{n(d^n+1)}
- h_{\crit,\bQ}\right|\\
\leq 8(d-1)(196d^2-192d-3)\cdot h^{\min}_{d,\bQ}\cdot\frac{\sum_{\ell|n}\sum_{m|\ell}m^2}{d^n+1}+C_{d,n}'
\quad\text{on }\rM_d(\overline{\bQ}).\tag{\ref{eq:local}$''$}\label{eq:multcrit}
\end{multline}

\subsection{Effective comparison among the ample, minimal, and critical height functions}

Choosing appropriate local one-to-finite multisections of 
the projection $\Rat_d\to\rM_d$,
Silverman's argument to compare the minimum/ample height functions
on $\rM_d(\overline{\bQ})$
based on the Weil height machine and on Siu's bigness criterion for the 
differences between divisors 
(for complex geometry, see e.g.\ the book
\cite{Lazarsfeld04}) is partly quantified
so that for some constant $A\in\bR$ depending 
only on $d$ and $\deg_D(\rM_d)$, 
the inequality between the minimum/ample height functions
\begin{gather}
h^{\min}_{d,\bQ}\le(2d-2)\cdot h_{\rM_d,D,\bQ}+A
\quad\text{on }\rM_d(\overline{\bQ})
\label{eq:minimalample}
\end{gather}
holds.
On the other hand, choosing similar (but more dynamical) 
kind of local one-to-finite multisections of 
the projection $\Rat_d\to\rM_d$,
by Theorem \ref{th:effectiveapprox} and Thurston's rigidity theorem,
there are infinitely many $n\in\bN$
such that the following comparison between multiplier/ample
height functions
\begin{gather}
 2C_1(d,D)\cdot h_{\rM_d,D,\bQ}-A_{d,n}
\leq \frac{h_{\mathbb{P}^{d^n+1},\bQ}\circ s_n}{n(d^n+1)}
\leq 2C_2(d,D)\cdot h_{\rM_d,D,\bQ}+A_{d,n}
\label{eq:multample}
\end{gather}
on $(\rM_d\setminus\rL_d)(\overline{\bQ})$ holds. Here
the constant $A_{d,n}\in\bR$ depends only on $d$ and $n$, and
the constants $C_1(d,D),C_2(d,D)>0$ are effectively computed
from the complex analytic quantities 
$\|\mu_{\bif,d}\|_{\rM_d(\bC)}$ and $\|T_{\bif,d}\|_{\rM_d(\bC)}$ and the algebraic $\deg_D(\rM_d)$.

\subsection{Concluding an effective solution of
Silverman's conjecture and an improvement of McMullen's finiteness theorem}
Now the effective comparisons/estimates 
\eqref{eq:multcrit}, \eqref{eq:minimalample}, 
and \eqref{eq:multample} concludes the comparison between $h_{\crit,\bQ}$
and $h_{\rM_d,D,\bQ}$ on $(\rM_d\setminus\rL_d)(\overline{\bQ})$
in Silverman's conjecture, having at least effective constants
$A_1,A_2>0$ in the conjecture (see Subsection \ref{th:silvermanconj}). 

This effective version of Silverman's conjecture 
together with \eqref{eq:multcritexact} and \eqref{eq:minimalample}
concludes the above 
improvement of McMullen's finiteness theorem (Theorem \ref{th:finiteness})
by an argument involving Northcott's finiteness theorem from arithmetic and a standard argument on field extensions.

\section{Meromorphic and hybrid families of dynamics
of rational functions}\label{sec:degenerate}

In this final section, let us see the case where
non-archimedean dynamics on Berkovich spaces
apply to the study of degeneration of a family of 
complex or non-archimedean dynamics; such a methodology at least
goes back to the study of
character varieties (see Morgan--Shalen \cite{MS85}).

\subsection{Meromorphic and hybrid families} 
Let $K$ be an algebraically closed field
of characteristic $0$ that is complete with
respect to a non-trivial norm $|\cdot|$.
The ring of $(K$-)analytic functions on
$\bD_K=\{t\in K:|t|<1\}$ is denoted by $\cO(\bD_K)$,
and we denote\footnote{following an idiomatic expression in relevant literature}
by $\cO(\bD_K)[t^{-1}]$
the ring of meromorphic functions on $\bD_K$
having no poles on $\bD_K^*:=\bD_K\setminus\{0\}$. An element $f\in(\cO(\bD_K)[t^{-1}])(z)$
of degree $d$ can be denoted by $(f_t)_{t\in\bD_K^*}$, 
where each $f_t(z)\in\bC(z)$ is the specialization of $f$ at $t\in\bD_K^*$,
and is called a meromorphic family of
rational functions on $\bP^1(K)$ parametrized by $\bD_K$
if in addition 
$\deg(f_t)=d$ for every $t\in\bD_K^*$. Then
we say $f$ degenerates at $t=0$ if $\deg(f_0)<d$,
where
writing $f=P(z)/Q(z)$ over the ring $\cO(\bD_K)[t^{-1}]$,
the rational function $f_0(z)$ on $\bP^1(\bC)$
of degree $\le d$ is the reduced ratio between
the specializations $P_0(z)$ and $Q_0(z)$ of $P,Q$ at $t=0$.

First, Theorem \ref{th:effectiveapprox} concludes
an asymptotic of $L(f_t)$ as $t\to 0$, which is
due to DeMarco \cite{DeMarco16} for archimedean
$K\cong\bC$, as follows.

\begin{mainth}[Gauthier--Okuyama--Vigny {\cite{GOV18}}]\label{th:degenerate}
Let $K$ be an algebraically closed field
of characteristic $0$ that is complete with
respect to a non-trivial norm $|\cdot|$.
Fixing an integer $d>1$,
for any meromorphic family $(f_t)_{t\in\bD_K^*}$
of rational functions on $\bP^1(K)$
of degree $d>1$ parametrized by $\bD_K$,
there is a constant $\alpha\ge 0$ such that
\begin{gather*}
 L(f_t)
=\alpha\cdot\log|t^{-1}|+o(\log|t^{-1}|)\quad\text{as }t\to 0.
\end{gather*}
\end{mainth}

Let $\bL$ denote the Levi--Chivita field,
which is the completion of an algebraic closure $\overline{\bC((t))}$
of the field $\bC((t))$ of formal Laurent series equipped
with a $t$-adic norm $|\cdot|_r$, $r\in(0,1)$, normalized as
$|t|_r=r$). Writing $\bD=\bD_{\bC}$
from now on, we regard a meromorphic family
$f\in(\cO(\bD)[t^{-1}])(z)$ of (complex) rational functions
on $\bP^1(\bC)$ of degree $d>1$ parametrized
by $\bD$ as an element
of $\Rat_d(\bL)$. In the case that
$f$ degenerates at $t=0$, the existence of
the weak limit
\begin{gather*}
 \lim_{t\to 0}\mu_{f_t}=:\mu_0\quad\text{on }\bP^1(\bC)
\end{gather*}
as well as various properties of $\mu_0$ including
the identification of $\mu_0$ with the pushforward under the reduction 
projection $\sP^1(\bL)\to\bP^1(\bC)$ (modulo $\mathfrak{m}_{\bL}$) of the
$f$-equilibrium measure say $\nu_f$ (rather than say $\mu_f$) on $\sP^1(\bL)$
is established by DeMarco \cite{DeMarco05} and 
DeMarco--Faber \cite{DF14}\footnote{for a complementation, see \cite{Okudegenerating}} 
(see Ma\~n\'e \cite{Mane88}
in the non-degenerating case, where
$\mu_0=\mu_{f_0}$). A more insight on the degenerating limit $\mu_0$ on $\bP^1(\bC)$
would be obtained by Favre \cite{Favre16} in terms of the hybrid space
(binding $\bP^1(\bC)\times\bD^*$ and $\sP^1(\bC((t)))$
into a new Berkovich space)
introduced by Berkovich and further developed by
Boucksom--Jonsson \cite{BJ17} (see also
Odaka \cite{Odaka19} etc.), where
the hybrid family of dynamical systems
is obtained from the family $f$ by replacing
the possibly degenerating dynamical system 
$(f_0,\bP^1(\bC))$ with the non-archimedean
dynamical system $(f,\sP^1(\bC((t))))$. Moreover,
Favre \cite{Favre16} also established the
asymptotic of $L(f_t)$ similar to that in Theorem
\ref{th:degenerate}, which even asserts that 
for any meromorphic family
$f$ of endomorphisms of $\bP^N(\bC)$, $N\ge 1$,
there is a constant $\alpha\ge 0$ such that
\begin{gather*}
 L(f_t):=\int_{\bP^N(\bC)}\log|\det Df_t|\mu_{f_t}
=\alpha\cdot\log|t^{-1}|+o(\log|t^{-1}|)
\quad\text{as }t\to 0,
\end{gather*}
and identified the constant $\alpha$ 
as the non-archimedean
Lyapunov exponent $L_{\operatorname{NA}}(f)$ of $f$
with respect to the $f$-equilibrium measure $\nu_f$ on $\sP^N_{\bC((t))}$ 
(up to the choice of $r$ for the $t$-adic norm $|\cdot|_r$ on $\bC((t))$). 

The error term $o(\log|t^{-1}|)$ in the above
asymptotic of $L(f_t)$ as $t\to 0$
is a $\rd\rd^c$-potential of the bifurcation measure on $\bD^*$,
which is a continuous function $\bD^*$, 
for the holomorphic family $(f_t)_{t\in\bD^*}$
of rational functions of degree $d$ on $\bP^1(\bC)$ 
(cf.\ Subsection \ref{sec:Lyap}), and
extends at least subharmonically to $\bD$ (so is bounded from above
around $t=0$). A further understanding of the asymptotic of 
this error term as $t\to 0$ is desirable, e.g.,\,in the study of geometries of the dynamical moduli $\rM_d$. 

In the $N=1$
dimensional (and even non-archimedean $K$)
case, Favre--Gauthier \cite{FG18}
established the ($\bR$-valued) continuous
extension across $t=0$
of the error term $o(\log|t^{-1}|)$
for any meromorphic family $f\in(\cO(\bD_K)[t^{-1}])[z]$ 
of polynomials of degree $d>1$.
In general, the situation is more complicated.

\begin{mainth}[DeMarco--Okuyama \cite{DO17}]\label{th:recipe}
There is a recipe to construct
a degenerating meromorphic family $(f_t)_{t\in\bD^*}$
of rational functions on $\bP^1(\bC)$ of degree $d>1$ 
parametrized by $\bD$ such that
\begin{gather*}
 \lim_{t\to 0}\bigl(L(f_t)-\alpha\cdot\log|t^{-1}|\bigr)=-\infty.
\end{gather*}
\end{mainth}
We constructed several kinds of such examples $f$
using the recipe in Theorem \ref{th:recipe}.

\subsection{Degeneration and bifurcation}

In the case of $N=2$, 
the space $\operatorname{Hol}_2(\bC\bP^2)$
of quadratic holomorphic
endomorphism of $\bC\bP^2$ is identified with
a complex hypersurface complement of $\bC\bP^{17}$
($17=3\cdot 4!/(2!2!)-1$) by coefficient parametrization,
and all (normalized) quadratic H\'enon maps
\begin{gather*}
 (z,w)\mapsto(w,cz+w^2+c_1w+c_2),\quad c\in\bC^*,c_1,c_2\in\bC,
\end{gather*}
which are the only dynamically non-trivial
quadratic polynomial automorphisms of $\bC^2$,
live in the complement hypersurface in $\bC\bP^{17}$.

\begin{mainth}[Bianchi--Okuyama \cite{BO18}]
The H\'enon map locus in the hypersurface
$\bC\bP^{17}\setminus\operatorname{Hol}_2(\bC\bP^2)$
is contained in the closure in $\bC\bP^{17}$
of the $J$-unstable (bifurcation) locus in
$\operatorname{Hol}_2(\bC\bP^2)$.
\end{mainth}

One of the key ingredients in the proof is an argument similar to that
in the proof of the (not only continuous but also) harmonic extension
across $t=0$ of the error term $L(f_t)-(1/2)\log|t^{-1}|$ (also having
$\alpha=1/2$ in this setting)
for the meromorphic family of quadratic holomorphic endomorphisms
\begin{gather*}
 f_t(z,w)=f_t(z,w;g,h)=\begin{pmatrix}
	   w\\
	   cz+w^2+c_1w+c_2
	  \end{pmatrix}+t\begin{pmatrix}
			  g(z,w)\\
			  h(z,w)
			 \end{pmatrix},\quad 0<|t|\ll 1,
\end{gather*}
of $\bC\bP^2$ degenerating to the H\'enon map $(z,w)\mapsto(w,cz+w^2+c_1w+c_2)$
as $t\to 0$, where $(g,h)\in(\bC[z,w])^2$ satisfying
$\deg g=2, g_{zz}\in\bC^*$, and $\deg h\le 2$,
under the non-exceptionality assumption $h_{zz}/g_{zz}\neq c$.

\begin{acknowledgement}
The author thanks the referee and the editors
for careful scrutiny and helpful comments
on arguments and presentation in the drafts
of the original Japanese version of this expository article. 
The notes prepared in each occasion of
the author's giving talks in seminars, colloquiums,
conferences, minicourses, etc.,\ were very helpful
in preparing this expository article, and the author also 
thanks the organizers of all those occasions. In translation,
the author also thanks Professors Yohsuke Matsuzawa, Kaoru Sano, Thomas Gauthier,
and Gabriel Vigny for their comments.
\end{acknowledgement}


\def\cprime{$'$}

\end{document}